%% file: convex.09.tex
\documentclass[12pt]{article}
\usepackage{amsmath,amssymb,amsthm,amsfonts,graphicx,color,enumerate,float,epsfig}
\begin{document}

\newtheorem{ttt}{TTT}
\newtheorem{thm}{Theorem}
\newtheorem{lemma}[thm]{Lemma}
\newtheorem{cor}[thm]{Corollary}
\newtheorem{prop}[thm]{Proposition}
\newtheorem{exercise}[ttt]{Exercise}
\newtheorem{question}[ttt]{Question}
\newtheorem*{mainthm0}{Theorem \ref{free'}}
\newtheorem*{mainthm}{Theorem \ref{mod}}
\newtheorem*{mainthm2}{Theorem \ref{mod2}}

\theoremstyle{remark}
\newtheorem{remark}[thm]{Remark}
\newtheorem{important remark}[thm]{Important Remark}
\newtheorem{definition}[thm]{Definition}
\newtheorem{q}[thm]{}
\newtheorem{example}[thm]{Example}
\newtheorem{fact}[thm]{Fact}
\newtheorem{convention}[thm]{Convention}

\def\diam{\operatorname{diam}}
\def\cal{\mathcal}
\def\R{{\mathbb R}}
\def\Z{{\mathbb Z}}

\def\H{{\mathbb H}}
\def\Q{{\mathbb Q}}
\def\F{{\mathbb F}}
\def\E{{\Bbb E}}
\def\tr{{\rm tr}}
\def\dcap{\cap_{10\delta}}
\def\Mod{{\rm Mod}}
\def\ep{{\varepsilon}}
\def\Isom{{\rm Isom}}
\def\C{{\mathcal C}}

\def\red{\textcolor{red}}
\title{Subgroups generated by two pseudo-Anosov elements
 in a mapping class group. II. \\
Uniform bound on exponents}
\author{Koji Fujiwara \\
Graduate School of Information Science \\
Tohoku University \\
Sendai, 980-8579, Japan \\
{\tt fujiwara@math.is.tohoku.ac.jp}
}


\maketitle




\abstract{Let $S$ be a compact orientable surface, and
$\Mod(S)$ its mapping class group.
Then there exists a constant $M(S)$, which 
depends on $S$,  with the following property.
Suppose $a,b \in \Mod(S)$ are independent 
(i.e., $[a^n,b^m]\not=1$ for any $n,m \not=0$) pseudo-Anosov elements.
Then for any $n,m \ge M$, the subgroup
$\langle a^n,b^m \rangle$ is free of rank two, and 
convex-cocompact in the sense of Farb-Mosher. In particular 
all non-trivial elements in $\langle a^n,b^m \rangle $
are pseudo-Anosov.
We also show that there exists a constant $N$, which depends on 
$a,b$, such that 
$\langle a^n,b^m \rangle$ is free of rank two 
and convex-cocompact if $|n|+|m| \ge N$ and $nm \not=0$. 
}


\section{Introduction}
This is the second half of our study
on subgroups generated by two pseudo-Anosov 
elements in a mapping class group.
We will improve the results we obtained in the first
half of the study \cite{Fu}, but 
one can read this paper independently.
We explain the improvement after we state
the main results in this section.

\subsection{Hyperbolic isometry and (quasi-)axis}
A geodesic space is called $\delta$-{\it hyperbolic} 
for $\delta \ge 0$ if for any geodesics $\alpha,\beta,\gamma$
which form a triangle, $\alpha$ is contained
in the $\delta$-neighborhood of $\beta \cup \gamma$
(\cite{Gr.hyp}).
Let $\Gamma$ be a $\delta$-hyperbolic graph.
Let $a$ be an isometry of $\Gamma$.
If there exist a point $x \in \Gamma$ and a constant $C >0$ such that 
$d(x,a^n(x)) \ge Cn$ for any $n>0$, then $a$ is 
called {\it hyperbolic}.

Suppose $a$ is a hyperbolic isometry.
If there exists a geodesic $\alpha$ 
such that $a(\alpha)$ is contained in the $C$-neighborhood
of $\alpha$ for some $C \ge 0$,
$\alpha$ is called a {\it quasi-axis} of $a$.
By $\delta$-hyperbolicity of $\Gamma$, it then follows
that $a(\alpha)$ is in the $2\delta$-neighborhood of $\alpha$.
If $\alpha$ and $\beta$ are quasi-axes of $a$ (they are 
geodesics by definition), then they are contained
in the $2\delta$-neighborhood of each other.
If $C=0$ we say $\alpha$ is an {\it axis}.
We remark that an axis and even a  quasi-axis may not exist for $a$,
but there is always a {\it quasi-geodesic} 
which is invariant by $a$. That will be  good enough for our argument
(See section \ref{no.axes}).
In the literature, quasi-axis sometimes means
a quasi-geodesic which is invariant by $a$.
Our definition is different.

For two points $x,y \in \Gamma$, we may 
denote a geodesic joining them by $[x,y]$ although
they are not unique.
We may write the distance between the two points
by $|x-y|$.

For an isometry $a$, we define its {\it translation length},
$\tr(a)$, by 
$$\tr(a) =\lim_{n \to \infty} \frac{|x-a^n(x)|}{n} \ge 0$$
for a point $x$. It is easy to see $\tr(a)$ does not
depend on the choice of $x$.
The isometry $a$ is hyperbolic iff $\tr(a)>0$.

It is known (\S 7,8 \cite{Gr.hyp})
that $a$ is hyperbolic 
if there is a point $p \in \Gamma$ such that
the following is satisfied.
In particular, the element $a$ has infinite order.
We remark that if $\delta >0$ we can replace $(\delta +1)$ by $\delta$.
\begin{equation}\label{hyp}
|p-a(p)| \le |a(p)-a^{-1}(p)|-100(\delta+1).
\end{equation}
The argument is geometric. 
Consider the following infinite path,
which is $a$-invariant.
$$ \beta= \cdots [a^{-1}(p),p] \cup   [p,a(p)] 
\cup [a(p),a^2(p)] \cup [a^2(p),a^3(p)] \cup \cdots $$
For each $n>0$, by $\delta$-hyperbolicity,
$[p,a^n(p)]$ is contained in the $n \delta$-neighborhood
of $[p,a(p)] \cup \cdots \cup [a^n(p),a^{n+1}(p)]$.
Using (\ref{hyp}), one can show 
it is indeed contained in the $3\delta$-neighborhood.
It then follows that $\beta$ is a quasi-geodesic.
If $a$ has a (geodesic) quasi-axis, then it is 
contained in the $4\delta$-neighborhood of $\beta$.

In section \ref{section.freesubgroup} we use a similar geometric 
idea (Proposition \ref{gromov.qi}) to give a sufficient condition
for two isometries to generate a free group.

\subsection{Convex co-compact subgroup in $\Mod(S)$}
Suppose that a finitely generated group $G$ is acting 
on $\Gamma$ by isometries.
Fix a finite generating set and let $|a|$
be the word metric of $a \in G$.
Let $x \in \Gamma$ be a point and consider
the map from $G$ to $\Gamma$ defined by sending 
$a \in G$ to $a(x) \in \Gamma$.
We call this map as {\it the embedding by an orbit}
of the action by $G$.
If there exist constants $L,C>0$ such that 
for all $a \in G$
$$|a|/L-C \le d(x,a(x)) \le L|a|+C,$$
then we say the map is {\it quasi-isometric}.

Our main application is regarding subgroups
of mapping class groups.
Let $S$ be a compact orientable surface, 
and $\Mod(S)$ its mapping class group.
Let $\C(S)$ be the curve graph of $S$, on which 
$\Mod(S)$ acts by isometries
(see for example \cite{Iv}, \cite{MM} for the definition).

Masur-Minsky \cite{MM} showed that
$\C(S)$ is $\delta$-hyperbolic and 
an element $a \in \Mod(S)$
is pseudo-Anosov if and only if it
acts as a hyperbolic isometry on $\C(S)$, and
moreover that there always exists a quasi-axis.

For a subgroup $G < \Mod(S)$, Farb-Mosher \cite{farb.mosher}
introduced the notion of {\it convex-cocompact}
in terms of the action on Teichmuller space.
It has been shown (\cite{Ha}, \cite{KeLe}) that $G$ is convex-cocompact iff
for a point $c \in \C(S)$, the map from $G$ to $\C(S)$
sending $g$ to $g(c)$ is quasi-isometric.
Note that the choice of the generating set 
and the point $c$ is not important.

\subsection{Main results}
In Section \ref{section.freesubgroup}
we discuss subgroups generated
by powers of two hyperbolic isometries on a hyperbolic graph,
and obtain sufficient conditions for them to be free.
That section is the main technical  part of the paper.
We put an overview of the argument
in Section \ref{overview}.
The final result is the following.
The point is that although the constant $M$ depends
on $\Gamma$ and the action of $G$ on $\Gamma$, it does not depend
on $a$ and $b$. It will become clear how the constant
$M$ depends on the action.

\begin{mainthm0}
Suppose $G$ acts acylindrically on a $\delta$-hyperbolic graph $\Gamma$.
Then there exists a constant $M$
with the following property.

Suppose $a,b\in G$ act hyperbolically.
Assume for any $p,q \not=0$,
$[a^p,b^q] \not= 1$ in $G$.
Then for any $n,m \ge M$,
$\langle a^n,b^m \rangle$ is free of rank two.
Moreover, the embedding of $\langle a^n,b^m \rangle$ by an orbit 
in $\Gamma$ is quasi-isometric.
In particular, all non-trivial elements
in $\langle a^n,b^m \rangle$ are hyperbolic on $\Gamma$.
\end{mainthm0}

The {\it acylindricity} of an action (see section \ref{section.acylindrical}
for the definition) is a weak assumption
on properness, and in particular, the result applies
to a word-hyperbolic group and its action on a Cayley graph,
therefore, if $G$ is a word-hyperbolic group, then there
exists $M$ such that for any two elements $a,b \in G$ of 
infinite order, either the subgroup
$\langle a, b \rangle$ is elementary or else
for any $n,m \ge M$, $\langle a^n, b^m \rangle$
is free and quasi-convex in $G$.
It seems this claim is new
(see \cite[8.2 E]{Gr.hyp} for the statement
without a bound on $n,m$).
The result also immediately applies to mapping class groups. 
That is our motivation and we show the following 
in Section \ref{section.mcg}.

\begin{mainthm}
Let $S$ be a compact orientable surface, and
$\Mod(S)$ its mapping class group.
Then there exists a constant $M(S)$ with the following property.
Suppose $a,b \in \Mod(S)$ are pseudo-Anosov elements
such that $[a^n,b^m]\not=1$ for any $n,m \not=0$.
($a,b$ are called independent.)
Then for any $n,m \ge M$, the subgroup 
$\langle a^n,b^m \rangle$ is free of rank two, and 
convex-cocompact in the sense of Farb-Mosher. In particular 
all non-trivial elements in $\langle a^n,b^m \rangle $
are pseudo-Anosov.
\end{mainthm}

It was known (\cite{Iv}, \cite{Mc}) that  $\langle a^n,b^m \rangle$
is free for sufficiently large $n,m$.
A uniform bound on $n,m$ is new.
In our previous study \cite{Fu}, a uniform bound on one of $n,m$ was shown.
Namely, there exists a constant $N(S)$ such that 
 $\langle a^n,b^m \rangle$ is free and convex-cocompact if one of $n,m$ is 
at least $N$ and the other one is sufficiently large.
The previous study has also been 
used to show the {\em uniform exponential growth}
of a mapping class group by Mangahas \cite{Man}.

The theorem concerns only pseudo-Anosov elements $a,b$.
It is unknown if a uniform bound such as  $M(S)$
exists for two elements $a,b$ of infinite order in general
such that  the subgroup $\langle a^n,b^m \rangle$ is free if $n,m \ge M$.
Note that the subgroup is never convex-cocompact
unless both $a$ and $b$ are pseudo-Anosov.
In the case that both $a,b$ are Dehn twists (\cite{Is}), and more generally,
positive multi-twists (\cite[Theorem 3.2]{Ham}),
it is known that  $\langle a^n,b^m \rangle$ is 
free or abelian if $n,m \ge 2$.
Recently, Leininger and Margalit \cite{LeMa}
have shown that if $S$ is the $n$-times punctured sphere,
then for any two elements $a,b \in \Mod(S)$,
$\langle a^N,b^N \rangle$ is either free or abelian
for $N=n!$.

It would be interesting to know 
for which $(n,m)$, the subgroup $\langle a^n,b^m \rangle$
is free for given $a,b$ in the above theorem.
The following theorem says that for given $a,b$, 
the subgroup $\langle a^n,b^m \rangle$
is free except for finitely many pairs $(n,m)$.
It is not clear if the number of those exceptional pairs
$(n,m)$ is bounded, but we know that 
the constant $N$ depends on $a,b$ 
in the following theorem (see Example \ref{example}).
\begin{mainthm2}
Let $S$ be a compact orientable 
surface and $a,b$ two independent
pseudo-Anosov elements.
Then there exists $N$ such that for any
$n \ge N$, both 
$\langle a,b^n \rangle$ and  $\langle b,a^n \rangle$
are free of rank two, and convex-cocompact.
In particular, 
$\langle a^n,b^m \rangle$ is 
free of rank two, and convex-cocompact
if $|n|+|m| \ge 2N$ and $nm\not=0$.
\end{mainthm2}


The author would like to thank Z. Sela for many insightful
suggestions, and L. Mosher for his interest and comments.
He is grateful to M. Bestvina and T. Delzant.
The work was partly done during his stay at MSRI
in Fall 2007. He appreciates their hospitality.
He wishes to acknowledge the support of
JSPS-KAKENHI-19340013.

\section{Free subgroups}\label{section.freesubgroup}
\subsection{Overview}\label{overview}
In this section, we will find sufficient
conditions for certain powers of two {\it independent}
hyperbolic isometries $a,b$ to generate a free group.
The final results are Theorem \ref{free} and Theorem \ref{free'}.

It is well-known that for sufficiently large $n,m>0$,
$a^n,b^m$ generate a free group (Proposition \ref{nielsen}).
The argument is an application of a geometric fact
on a $\delta$-hyperbolic space 
(Proposition \ref{gromov.qi}).
The goal of this section is to give an upper bound
on $n,m$ which does not depend on $a$ and $b$.

In Section \ref{section.acylindrical},
by analyzing the argument for Proposition \ref{nielsen} carefully,
we first show that there is 
an upper bound on $n$ and $m$
if the translation length of $a$ and $b$
are comparable (Proposition \ref{nielsen1}).
A more difficult case is that 
one of the translation length, say for $b$,
is much smaller than the translation length 
of $a$. In this case, if we use the 
same argument, we need to take the exponent $m$
for $b$ very large so that $a^n,b^m$
generate a free group.
In Section \ref{section.nielsen2}, we use a different
idea to deal with this case and show
there is an upper bound on $n$ and $m$
if the translation of one of $a$ and $b$ is much smaller
than the translation length of the other
(Proposition \ref{nielsen2}).
This part is tedious (the idea is elementary,
but we put many details), but we think 
this is a main technical achievement of the paper. 

As we explained, the two 
propositions are complimentary to each other, and 
combining them, we obtain an upper bound
on the both exponents which does not depend on 
$a$ and $b$ in Section \ref{section.upperbound}
(Theorem \ref{free}).

We first prove those results  under the assumption 
that $a$ and $b$ have {\it quasi-axes}. Then in section 
\ref{no.axes} we explain that 
the assumption is indeed redundant, which 
gives Theorem \ref{free'}.
For our application in this paper, 
we only need Theorem \ref{free}, but
we prove Theorem \ref{free'} for potential 
application in the future. 

Another technical issue is the properness of an 
action. We argue under the assumption 
of {\it acylindricity}, which is weaker than 
the action being properly discontinuous
(see Section \ref{section.acylindrical}).
We will need that when we discuss application 
in Section \ref{section.mcg}.


\subsection{Nielsen condition}
In this section, we review 
a well-known
fact (Proposition \ref{nielsen}) and its 
proof. We start with 
a fundamental result from \cite{Gr.hyp} (see \S 7 and 8).

%

\begin{prop}[7.2C \cite{Gr.hyp}, Three points condition]\label{gromov.qi}
Let $\Gamma$ be a $\delta$-hyperbolic graph. 
Let $\ep > 100 \delta$ be a constant.
Let $p_i \in \Gamma (i \ge 1)$ be points 
such that for all $i \ge 1$
\begin{equation}\label{three.points}
|p_i-p_{i+2}| \ge \max (|p_i-p_{i+1}|,|p_{i+1}-p_{i+2}|) + \ep.
\end{equation}
Then, for each $i \ge 3$,
\begin{equation}\label{additive}
\lambda |p_1-p_i| \ge \sum_{j=1}^{i-1}|p_j-p_{j+1}|,
\end{equation}
where 
$\lambda=(\frac{\ep}{100}-\delta)^{-1} \max_{1 \le j \le i-1}|p_j-p_{j+1}|$.
In particular, $|p_1-p_i|>0$ for all $i \ge 3$.
\end{prop}

We call the inequality (\ref{three.points}) {\it the three points
condition} (for $p_i,p_{i+1},p_{i+2}$ and a constant
$\varepsilon$).
If the three points condition is satisfied
for any three consecutive points in a sequence,
then we say the sequence satisfies the three point
condition.


Proposition \ref{gromov.qi} has been used to derive a condition 
for sufficiently large powers of 
hyperbolic isometries $a,b$ of $\Gamma$ with quasi-axes 
 $\alpha, \beta$ to generate
a free group in terms of $\alpha, \beta$
(Proposition \ref{nielsen}). 
In that argument, it will be important 
how much of $\beta$ is contained in the 
$10 \delta$-neighborhood of $\alpha$ and vise-versa.
We thus define {\it the $10 \delta$-overlap} of $\alpha$ and $\beta$,
denoted by $\alpha \dcap \beta$, as follows.
$$\alpha \dcap \beta = (\alpha \cap N_{10 \delta}(\beta))
\cup (\beta \cap N_{10 \delta}(\alpha)).$$
Let $|\alpha \dcap \beta|$ denote the diameter of this set.
$|\alpha \dcap \beta|$ can be $\infty$.
If it is finite, by the $\delta$-hyperbolicity, the longest segment
of $\alpha$, the longest segment of $\beta$
and the longest geodesics which are contained in $\alpha \dcap \beta$
all have length between $|\alpha \dcap \beta|-20 \delta$
and $|\alpha \dcap \beta|+20 \delta$,
and those segments are in the $20 \delta$-neighborhood
of each other.

The following fact is well-known (\cite{Gr.hyp}, see also \cite{Ko}).
Notice that the exponents $n$ and $m$ which satisfy
the inequalities depend on $a$ and $b$.
For readers who have not seen a proof, we give details
of the argument, since
we will generalize the statement using the same 
idea.

\begin{prop}[Nielsen condition]\label{nielsen}
Suppose isometries $a,b$ act hyperbolically on 
a $\delta$-hyperbolic graph $\Gamma$ with
quasi-axes $\alpha,\beta$, respectively.
Suppose 
$|\alpha \dcap \beta| < \infty$.
If $1 \le n,m \in \Bbb{Z}$ are such 
that 
$$\tr(a^n) \ge  |\alpha \dcap \beta| + 100 (\delta+1), \,
\tr(b^m) \ge  |\alpha \dcap \beta| + 100 (\delta+1) $$
then 
$\langle a^n,b^m \rangle < \Isom(\Gamma)$ is free of rank two.
Moreover, for a point $x \in \Gamma$,
the embedding of the subgroup $\langle a^n,b^m \rangle$ 
to $\Gamma$ by sending $w$ to $w(x)$ is quasi-isometric.
\end{prop}

\proof 
We use Proposition \ref{gromov.qi} with 
$\ep=100(\delta+1)$.
We first show that $\langle a^n,b^m \rangle$ is free,
then argue independently the moreover part. It turns out that 
for a certain choice of $x$, the embedding is not 
only quasi-isometric, but also bi-Lipschitz,
which implies that $\langle a^n,b^m \rangle$ is free.
In that sense, the first part is not necessary, but 
we hope it will make the whole argument more transparent
in this way. 

Set 
$$A=a^n, B=b^m.$$
We remark that $\alpha,\beta$ are quasi-axes of $A,B$, respectively.
Let $w$ be a non-empty reduced word on $A,B$ and 
we prove that the action of $w$ on $\Gamma$
is non-trivial, therefore, $w$ is a non-trivial element
in ${\rm Isom}(\Gamma)$. 
It suffices to find a point $p \in \Gamma$ with $w(p)\not=p$.
(In this proof the point $p$ does not depend on the word $w$).

\begin{figure}[htbp]
\begin{center}
\scalebox{0.7}{\input{fig0.pstex_t}}
\caption{Nielsen condition}
\label{figure0}
\end{center}
\end{figure}
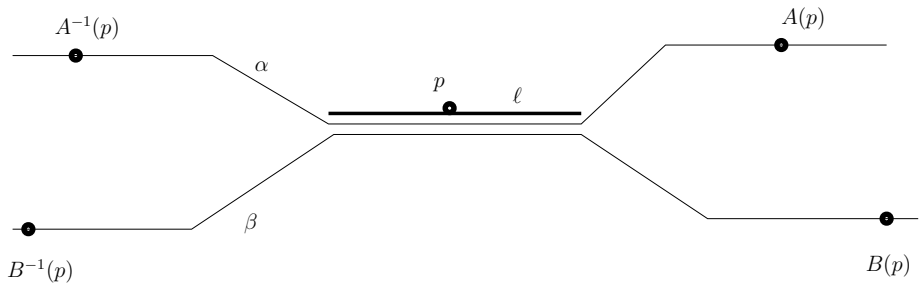

Suppose $|\alpha \dcap \beta| =0$.
Let $\ell$ be a geodesic segment which 
realizes the distance between $\alpha$
and $\beta$, and $p \in \ell$ 
the mid point. (See Figure \ref{figure0}.)
We claim $w(p)\not=p$.
To see it, let 
$$w=A^{n_1}B^{m_1} \cdots A^{n_i}B^{m_i},$$
where $n_1,m_i$ are possibly $0$.
We discuss the case such that both $n_1,m_i$ are not $0$.
Set
\begin{align*}
p_0&=p,
p_1=A^{n_1}(p),p_2=A^{n_1}B^{m_1}(p),
p_3=A^{n_1}B^{m_1}A^{n_2}(p), \cdots, \\
p_{2i-1}&=A^{n_1}B^{m_1}\cdots A^{n_i}(p), 
p_{2i}=A^{n_1}B^{m_1} \cdots A^{n_i}B^{m_i}(p).
\end{align*}
Then the the sequence of points $p_j$ satisfies the 
inequalities (\ref{three.points}) in Proposition \ref{gromov.qi}.
Indeed, for example, for $p_1,p_2,p_3$, if we apply 
the element $B^{-m_1}A^{-n_1}$, then we get $B^{-m_1}(p),p,A^{n_2}(p)$.
Since $m_1,n_2 \not=0$, 
those three points
satisfy the three point condition by the $\delta$-hyperbolicity
of $\Gamma$.
Now, by Proposition \ref{gromov.qi}, we get 
$|p_0-p_{2i}|>0$, namely $|p-w(p)|>0$.
We can argue similarly if $n_1$ or $m_i$ is $0$,
and omit the details.

Suppose  $|\alpha \dcap \beta| =D>0$.
Let $\ell$ be the longest geodesic 
which is contained in $\alpha \dcap \beta$.
Then $D-20 \delta \le |\ell| \le D+20 \delta$.
Let $p$ be the mid point of $\ell$.
Define the points $p_j$ in the same way as before, 
then by using the assumption on 
$\tr(a^n)$ and $\tr(b^m)$, 
we get $|p-w(p)|>0$ by Proposition \ref{gromov.qi}.

Now we argue that the embedding by the orbit of a point $x$ is quasi-isometric
with respect to the word metric for 
$a^n,b^m$. We remind that the choice of the point $x$ is not important.
Since $\langle a^n,b^m \rangle$ is finitely generated,
the embedding is always Lipschitz
with respect
to the constant $\max (|a^n(x)-x|,|b^m(x)-x|)$.
We will show that there exist a point $x$, 
constants $L>0,C\ge0$ such that for any non-trivial reduced word
$w$ on $a^n,b^m$,
$$ \frac{|w|}{L} -C \le |w(x)-x|,$$
where $|w|$ is the word metric with respect to $a^n,b^m$.
Indeed we will have $C=0$ if  $x=p$.

The argument is a modification 
of the previous one, so that we use Proposition \ref{gromov.qi}.
For each $j$ set
$$B_j=[p_{2j},p_{2j-1}], \, \, 
A_{j+1}=[p_{2j},p_{2j+1}].$$
They are geodesic segments. 
Note that 
$$|p_{2j}-p_{2j-1}|=|p-B^{m_j}(p)|=|m_j| \tr(B),
\, \, |p_{2j}-p_{2j+1}|=|p-A^{n_{j+1}}(p)|=|n_{j+1}| \tr(A).$$
Thus the right hand side of the inequality
(\ref{additive}) is proportional to
the word length of $w$ when we vary $w$.
The left hand side of the inequality
is $\lambda|p_0-p_{2i}|=\lambda|p-w(p)|$,
so that if there is an upper bound on $\lambda>0$ 
when we vary the word $w$, then we would be done.
But since $\lambda = \max_j|p_j-p_{j+1}|$, 
$\lambda$ can be arbitrarily large when we vary $w$.

As a remedy, we will divide each of 
geodesic segments $A_j,B_j$, namely,
introduce certain points on them 
such that there is an upper bound on the 
distance between any two consecutive points, 
then apply Proposition \ref{gromov.qi}
to this new sequence of points.
The point is that introducing new points
does not change the property that the right hand side
of the inequality (\ref{additive}) is proportional to 
$|w|$, while the constant $\lambda$ for the 
new sequence will have an upper bound.

We start the argument. 
Again, we discuss the case that both $n_1,m_i$ are not $0$.
First, between $p_0$ and $p_1$, if $n_1>0$, define points by 
$$p_{0,0}=p_0, p_{0,1}=A(p),p_{0,2}=A^2(p), \cdots,
p_{0,n_1}=A^{n_1}(p)=p_1,$$
and if $n_1 < 0$, define the points 
$p_{0,k}$ for $0 \le k \le -n_1$ similarly
using actions by $A^{-1},A^{-2}, \cdots, A^{-n_1}$ on $p$.
Next, between $p_1$ and $p_2$, if $m_1 >0$, define points by 
$$p_{1,0}=p_1, p_{1,1}=A^{n_1}B(p),p_{1,2}=A^{n_1}B^2(p), \cdots,
p_{1,m_1}=A^{n_1}B^{m_1}(p)=p_2,$$
and if $m_1 < 0$, define the points similarly as follows.
$$p_{1,0}=p_1, p_{1,1}=A^{n_1}B(p),p_{1,2}=A^{n_1}B^2(p), \cdots,
p_{1,|m_1|}=A^{n_1}B^{m_1}(p)=p_2.$$
We define points similarly between $p_{2j}$ and $p_{2j+1}$, 
and also between $p_{2j}$ and $p_{2j-1}$ for all $j$.
We obtain a sequence of points $p_{j,k}$ with the canonical 
order (the lexicographical order on $(j,k)$).
By definition, the distance between any two consecutive points
is either $\tr(A)$ or $\tr(B)$. 
Also, by our assumption, the sequence of points satisfies the 
three points condition 
of Proposition \ref{gromov.qi} such that 
$\lambda=\max(\tr(A),\tr(B))$, which no more depends on $w$.
Here, we regard, for example, $p_{0,n_1}$ and $p_{1,0}$ 
are the same point, which is $p_1$.
By the proposition (for the first inequality), we get 
\begin{align*}
|p_0-p_{2i}| & \ge \frac{1}{\lambda} \sum_{j=0}^{i-1}
\left(\sum_{k=1}^{|n_{j+1}|}|p_{2j,k-1}-p_{2j,k}|
+\sum_{k=1}^{|m_{j+1}|}|p_{2j+1,k-1}-p_{2j+1,k}|\right) \\
&= \sum_{j=0}^{2i-1}|p_{j}-p_{j+1}|
= \frac{1}{\lambda}\left(\tr(A)\sum_{j=1}^i|n_j| 
+ \tr(B)\sum_{j=1}^i|m_j|\right) \\
&\ge \frac{\min(\tr(A),\tr(B))}{\lambda}|w|.
\end{align*}
Set $L'=  \frac{\min(\tr(A),\tr(B))}{\lambda}>0$.
We get $|p-w(p)| \ge L'|w|$.
The same bound holds if $n_1$ or $m_i$ is $0$
as well.
We have shown that the embedding is bi-Lipschitz,
for this particular choice of a base point and 
a generating set.
\qed

We are interested in finding an upper bound
on $n,m>0$ such that $\langle a^n,b^m \rangle $ is not free
under some condition on the action of $G$
on $\Gamma$, provided that 
$D=|\alpha \dcap \beta| < \infty$.
Before we discuss that, 
we analyze the case when $D=\infty$,
namely, hyperbolic isometries $a$ and $b$ have a common quasi-axis.
For example, suppose $G$ is a word-hyperbolic group 
and $\Gamma$ is a Cayley graph, which is 
$\delta$-hyperbolic.
Then $D=\infty$ implies that 
$[a,b^k]=1$ for some $k>0$ and 
$[b,a^l]=1$ for some $l>0$,
where the commutator of two elements is defined by
$$[f,g]=f^{-1}g^{-1}fg.$$
To see the first claim, take a point $x \in \alpha$, the 
common quasi-axis, and 
look at the set of points
$$[a,b](x),[a,b^2](x),[a,b^3](x), \cdots$$
They are all in the $20 \delta$-neighborhood of $x$.
Since the action of $G$ is proper, there are only 
finitely many elements $g \in G$ with $|x-g(x)| \le 20 \delta$, 
therefore there must be distinct integers $n,m>0$
such that $[a,b^n]=[a,b^m]$.
This implies $[a,b^{n-m}]=1$.
The second claim is similar.

Notice that in the previous argument, we do not need that $D$ is infinite,
but it is enough if $D$ is sufficiently large.
To formulate a precise statement (Lemma \ref{axis}), we 
consider a certain condition, {\it acylindricity}, on the action
in the next section.

\subsection{Acylindrical action}\label{section.acylindrical}
In this section, we assume
certain properness of an action, acylindricity,
and improve Proposition \ref{nielsen}
to Proposition \ref{nielsen1}.

Let $\Gamma$ be a $\delta$-hyperbolic graph, 
and $G$ a group acting on $\Gamma$ by isometries.
Bowditch \cite{Bo}
defined that the action is {\it acylindrical} 
if for any $R>0$, there exist $K(R),L(R) \ge 1$ such that 
for any vertices $x,y \in \Gamma$ with $d(x,y) \ge L(R)$, 
the following set has at most $K(R)$ elements:
$$\{ g \in G| d(x,g(x)) \le R, d(y,g(y)) \le R \}.$$

\begin{lemma}\label{min.tr}
Suppose $G$ acts on a $\delta$-hyperbolic graph
$\Gamma$.
If the action is acylindrical with constants
$K(R), L(R)$, then there exists
an integer $P \ge 1$ such that for any element
$a \in G$ which acts hyperbolically on $\Gamma$
with a quasi-axis, we have $\tr(a^P) \ge 1$.
The constant $P$ depends only on $\delta$
and $K(200\delta)$.
\end{lemma}

\begin{convention}[Subscript of a constant]\label{convention}
To keep track of constants, we may number 
a constant by 
the number of the claim which the constant first appears, 
for example, the constant $P$ in 
Lemma \ref{min.tr} will be $P_{\ref{min.tr}}$.
We may omit the subscript if there is no confusion.
\end{convention}

\proof
If $\delta=0$, $\Gamma$ is a tree.
Then $\tr(a) \ge 1$. Set $P=1$.
Suppose $\delta >0$.
Set $R = 100 \delta$.
Let $\alpha$ be a (geodesic) quasi-axis of $a$.
Take a point $x \in \alpha$.
Let $y \in \alpha$ be a point with 
$|x-y| \ge L(2R)$.
If $|a^i(x)-x| \le R$ for some $i$, then 
$|a^i(y)-y)|\le 2R$.
This is because $a^i(\alpha)$ is in the $2\delta$-neighborhood
of $\alpha$.
By the acylindricity, this implies  that 
there is $I$ with $1 \le I \le K(2R)=K(200\delta)$
such that $|a^I(x)-x| >R$.
It then follows that for any $n \ge 1$, 
$|a^{In}(x)-x| > n(R-10\delta)$.
Thus 
$$\tr(a) \ge \frac{R-10\delta}{I} \ge \frac{90\delta}{K(200\delta)}.$$
Choose an integer $P$ such that 
$P \ge \frac{K(200\delta)}{90\delta}$.
\qed

\begin{lemma}\label{axis}
Let $\Gamma$ be a $\delta$-hyperbolic graph, 
and $G$ a group acting on $\Gamma$ acylindrically
with constants $K(R),L(R)$.
Suppose $a,b\in G$ act hyperbolically with quasi-axes 
$\alpha,\beta \subset \Gamma$, respectively. 

If  $a^n b\not= ba^n$ for all  $n \not=0$ or
$b^n a\not=a b^n$ for all  $n \not=0$, then 
$$|\alpha \dcap \beta| < 4P_{\ref{min.tr}}
K(20\delta) L(20\delta) \max(\tr(a),\tr(b))
+ 100 \delta.$$
\end{lemma}
By Convention \ref{convention},
$P_{\ref{min.tr}}$ is the  constant from Lemma \ref{min.tr}.

\proof
To argue by contradiction, suppose that 
the inequality was false.
Set $K=K(20 \delta), L=L(20 \delta)$.
For concreteness, suppose $\tr(b) \le \tr(a)$.
By our assumption, since $|\alpha \dcap \beta|$
is much larger than $2\delta$, 
the set $\alpha \dcap \beta$ looks like a narrow
tube. 

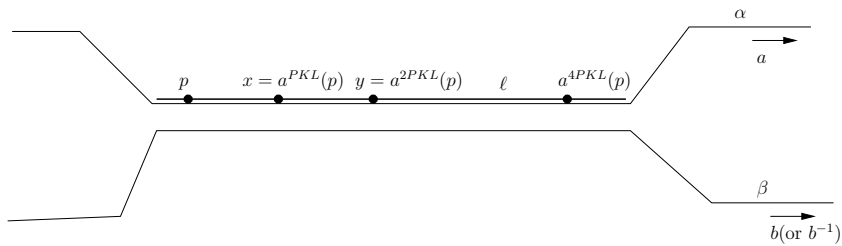
\begin{figure}[htbp]
\begin{center}
\scalebox{0.6}{\input{ueg.fig3.pstex_t}}
\caption{Apply the acylindricity to the pair $x,y$.}
\label{figure3}
\end{center}
\end{figure}

Let $\ell \subset \alpha$ be the longest segment
which is contained in $\alpha \dcap \beta$.
Then, by our assumption,  $|\ell| \ge 4PKL\tr(a) + 80 \delta$.
Take a point $p \in \ell$ such that 
the following points are in $\ell$ (See Figure \ref{figure3}.)
$$p, a(p), a^2(p), \cdots, a^{4PKL}(p).$$
Set
$$x=a^{PKL}(p), \, y=a^{2PKL}(p) \in \ell.$$
Since $y=a^{PKL}(x)$, and by Lemma \ref{min.tr}, 
$d(x,y) \ge PKL \tr(a) \ge KL \ge L$.

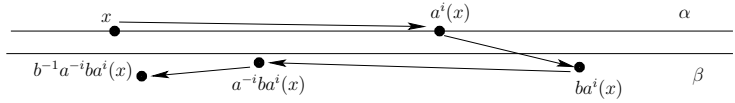
\begin{figure}[htbp]
\begin{center}
\scalebox{0.6}{\input{comm.pstex_t}}
\caption{How commutators $[b,a^i]$ act near $\alpha \dcap \beta$.}
\label{comm}
\end{center}
\end{figure}

\noindent
{\it Claim}. For each $i (1\le i \le PKL)$, 
$$d(x,[b,a^i](x)) \le 20 \delta, \, d(y,[b,a^i](y)) \le 20 \delta.$$
We first consider 
the special case that $\delta=0$, namely,
$\Gamma$ is a tree. Then, $\alpha \dcap \beta$ coincides 
the segment $\alpha \cap \beta$, and
also the segment $\ell$, therefore,  
all above points $a^n(p), 1 \le n \le 4PKL$, are in $\alpha \cap \beta$.
We want to show $x=[b,a^i](x)$, but 
this is obvious since when we apply
$a^{i},b,a^{-i}$, then $b^{-1}$ to $x$,
the point moves within $\ell$.
Thus, $[b,a^i](x)=x$.
If $ \delta >0$, we can show that the point moves 
in the $10 \delta$-neighborhood of $\ell$
when we apply $a^i,b,a^{-i}$ followed by $b^{-1}$ 
to $x$. Therefore, we get $d(x,[b,a^i](x)) \le 20 \delta$
by estimating the error terms
from the tree case using  triangle inequality. 
We leave the details to readers. (See Figure \ref{comm}.)
We can show $d(y,[b,a^i](y)) \le 20 \delta$
in the same way. We got the claim.

Since $|x-y| \ge L=L(20\delta)$, by the acylindricity of the action,
it follows from the claim that 
there are at most $K$ distinct elements
in the set $[b,a^i] (1 \le i \le PKL)$.
By the pigeon-hall principle, 
$[b,a^i]=[b,a^j]$ for some $i\not=j, (1 \le i,j \le PKL)$.
It follows that $a^ib^{-1}a^{-i}=a^jb^{-1}a^{-j}$,
therefore, $a^{i-j}b^{-1}=b^{-1}a^{i-j}$.
We get $[b,a^n]=1$ for some $n \not=0$.
The same argument applies to the elements $[a,b^i]$
since $\tr(b) \le \tr(a)$, therefore we also get
$[a,b^n]=1$ for some $n \not=0$ as well.
This is a contradiction.
\qed

Combining Proposition \ref{nielsen}
and Lemma \ref{axis}, we obtain the 
following.
This says that we can find a global bound
on one of the exponents, but the other bound
depends on the ratio $\tr(a)/\tr(b)$.

\begin{prop}\label{nielsen1}
Let $G$ be a group which acts on a $\delta$-hyperbolic 
graph $\Gamma$ acylindrically with constants $K(R),L(R)$.
Then, there exists a constant $N=N_{\ref{nielsen1}} \ge 1$
with the following property.
$N$ depends only on $\delta, K(20\delta),L(20\delta)$
and $K(200\delta)$.

Suppose $a,b \in G$ act hyperbolically
with quasi-axes.
Assume $[a^n,b^m] \not=1$ for all $n,m\not=0$.
Suppose there exists a number $q \ge 1$ such that 
$$\tr(a)/q \le \tr(b) \le \tr(a).$$
Then, for any $n \ge N$ and $m \ge qN$, 
$\langle a^n,b^m \rangle $ is free.
Moreover, the embedding of $\langle a^n,b^m \rangle $  
in $\Gamma$ is quasi-isometric.
\end{prop}

\proof
Put $K=K(20 \delta), L=L(20 \delta)$.
It then suffices to set 
$$N=4P_{\ref{min.tr}}KL+200(\delta+1)P_{\ref{min.tr}}.$$
Remember that $P_{\ref{min.tr}}$ depends
only on $\delta$ and $K(200\delta)$.
By Lemma \ref{axis}, since $\tr(b) \le \tr(a)$,
we have 
$|\alpha \dcap \beta| < 4PKL \tr(a) + 100 \delta$.
Since $\tr(a) \le q \tr(b)$, we have 
$|\alpha \dcap \beta| < 4qPKL \tr(b) + 100 \delta$.
Therefore, if $n \ge N$ and $m \ge qN$, 
then the conditions of 
Proposition \ref{nielsen}
are satisfied for $a^{n},b^{m}$, so that 
we conclude $\langle a^n,b^m \rangle$ is free, and 
the embedding is quasi-isometric.
\qed

\subsection{Another condition for freeness}\label{section.nielsen2}
In this section we discuss the case
when  $\tr(b)/\tr(a)$ is small
as opposed to Proposition \ref{nielsen1}.
We also apply Proposition \ref{gromov.qi} to 
a certain sequence of points, but we need 
a slightly different idea to construct the sequence
from a given word $w$.
The argument is 
elementary but lengthy, and takes
most part of this section.

\begin{prop}\label{nielsen2}
Suppose $G$ acts on a $\delta$-hyperbolic graph $\Gamma$
acylindrically with constants $K(R),L(R)$.
Assume that $f,g \in G$  act hyperbolically
with quasi-axes $\alpha,\beta$.
Assume that $[f^s,g^t]\not=1$ for any $s,t \not=0$.
Let $D=|\alpha \dcap \beta|$.
Suppose that the following conditions are satisfied.
\begin{enumerate}
\item
$\tr(g) \le \tr(f).$
\item
$D \le 2 \tr(f).$
\end{enumerate}
Then there exists a constant $N=N_{\ref{nielsen2}}>0$
with the following property. 
$N$ depends only on $\delta, K(20\delta),L(20\delta)$
and $L(200\delta)$.

If $n \ge N$, then 
$\langle g,f^n \rangle$ is free of rank two.
Moreover, the embedding 
of $\langle g,f^n \rangle$ in $\Gamma$ by an orbit is 
quasi-isometric.
\end{prop}

A few remarks are in order before the 
lengthy proof.
Unlike Proposition \ref{nielsen1}, 
the roles of two elements $f$ and $g$ are not symmetric.
For example, in the conclusion we take powers only for $f$.
Also, when we construct sequence of points
$p_i$ in the argument, we detect 
the action of $a=f^n$ more closely than the action of $b=g$, 
since $a$ moves a base point much more than $b$ does.
That is summarized as the condition (\ref{*}) in the claim.
As usual, in Part 1, we first show the subgroup $\langle g,f^n \rangle$
is free, by constructing a certain sequence of points $p_i$
which satisfies the three points condition 
because of the condition (\ref{*}). 
Then in Part 2, we show that the embedding of the subgroup by an orbit
is quasi-isometric by interpolating the points in the sequence $p_i$
as before. Part 2 is most complicated in the paper. 

\proof
First of all, $D < \infty$ by our assumption.
As in the proof for Proposition \ref{nielsen},
the argument is slightly different if $D=0$ than the case $D>0$, 
and from the view point of Proposition \ref{gromov.qi}, 
it is easier than the case $D>0$.
So, we discuss the case $D>0$ in detail, and then
discuss the case $D=0$ briefly.

Assume $D>0$. 

\noindent
{\bf Part 1}.
As usual, we first prove that $\langle g,f^n \rangle$ 
is free.
Set $K=K(20 \delta),L=L(20\delta)$.
This is the only place where constants
$K(R),L(R)$ are used.
Set 
$$E=D + 100 (\delta+1) + 10P_{\ref{min.tr}}KL \tr(f).$$
There exists a constant $N$, which depends only on 
$K,L,P,\delta$, such that if $n \ge N$, then 
$\tr(f^n) \ge 1000 E$. We used $\tr(f) \ge 1/P$
and $2 \tr(f) \ge D$.
This is the only place the condition 
$2 \tr(f) \ge D$ is used.
Fix such constant $N$.
Note that $N \ge 10000$
by the definition of $E$.

For $n \ge N$, set $a=f^n, b=g$.
We will show $\langle a,b \rangle $ is free.
Let $w$ be a non-empty reduced word on $a,b$,
and we show $w$ is non-trivial in ${\rm Isom}(\Gamma)$.
There are three cases according 
to the form of $w$.
In the case (O) it is clear that 
$w \not=1$ in ${\rm Isom}(\Gamma)$.
\begin{itemize}
\item[(O)] $w=b^m (m\not=0)$
\item[(I)]
$w=a^{n_1} b^{m_1} \cdots a^{n_i} b^{m_i} (i \ge 1)$
such that 
$n_1 \not=0, n_i \not=0$ and $m_i$ is possibly $0$,
\item[(II)]
$w=b^{m_0}a^{n_1} b^{m_1} \cdots a^{n_i} b^{m_i} (i \ge 1)$
such that $m_0 \not=0,  n_i \not=0$ and $m_i$ is possibly $0$.
\end{itemize}
For (I) and (II), we can find a point $p \in \Gamma $ such that $|p-w(p)|>0$, 
therefore $w$ is not $1$ in ${\rm Isom}(\Gamma)$.
The argument is very similar to each other, so
we only discuss the case (I) in detail.

Assume we are in the case (I).
Let $m$ be the mid point of a longest segment $\ell$
which is contained in $\alpha \dcap \beta$.
We can take a point $x \in \alpha$ and a point $y \in \beta$
such that $|x-y| \le 2 \delta$ and $|m-x|, |m-y| \le 4 \delta$.
It would be easier to follow the discussion 
if we imagine that $\alpha$ and $\beta$ coincide
in $\ell$ and that $m=x=y$, although 
there are actually errors of the order of $\delta$.

To show $|x-w(x)|>0$, 
we interpolate $x$ and $w(x)$ by 
the following points, which gives a sequence
satisfying the three points condition.

\begin{align*}
s_0 &=y
\\
p_1 &=x,\,q_1=a^{n_1}(x), \, r_1=a^{n_1}(y), \, s_1=a^{n_1}b^{m_1}(y), 
\\
p_2&= a^{n_1}b^{m_1}(x), \, q_2=a^{n_1}b^{m_1}a^{n_2} (x),
r_2=a^{n_1}b^{m_1}a^{n_2} (y), \, s_2=a^{n_1}b^{m_1}a^{n_2} b^{m_2}(y)
\\
& \cdots
\\
p_i&=a^{n_1}b^{m_1} \cdots a^{n_{i-1}} b^{m_{i-1}} (x), \, 
q_i=a^{n_1}b^{m_1} \cdots a^{n_{i-1}} b^{m_{i-1}} a^{n_i}(x),
\\
r_i&=a^{n_1}b^{m_1} \cdots a^{n_{i-1}} b^{m_{i-1}} a^{n_i}(y), \,
s_i=a^{n_1}b^{m_1} \cdots a^{n_{i-1}} b^{m_{i-1}} a^{n_i} b^{m_i}(y),
\\
p_{i+1}&= a^{n_{i-1}} b^{m_{i-1}} a^{n_i} b^{m_i}(x).
\end{align*}

Now set for each $i \ge j \ge 1$
$$A_j=[p_j,q_j],B_j=[r_j,s_j].$$

\noindent
{\bf Claim}. The following conditions
are satisfied for all $1 \le j$.

\begin{equation}\label{*}\tag{*}
\begin{cases}
(1). \, |q_j-r_j|\le 2 \delta, |p_j-s_{j-1}| \le 2 \delta.
\\
(2). \, |A_j| \ge 1000 E.
\\
(3). \, |A_j \dcap B_j| \le E. 
\\
(4). \, |B_{j-1} \dcap A_{j}| \le E. 
\\
(5). \, |A_j \dcap A_{j+1}| \le E.
\end{cases}
\end{equation}

We verify those later, and proceed to show the set of the 
conditions (1)--(5) implies 
that the sequence $\{p_j\}$ satisfy the three points 
condition for $990E$, therefore we get  $|p_1-p_{i+1}| >0 $, namely, 
$|w(x)-x| >0$. 

\begin{figure}[htbp]
\begin{center}
\scalebox{0.8}{\input{fig1.pstex_t}}
\caption{$A_j \dcap A_{j+1} \not= \emptyset$}
\label{figure1}
\end{center}
\end{figure}

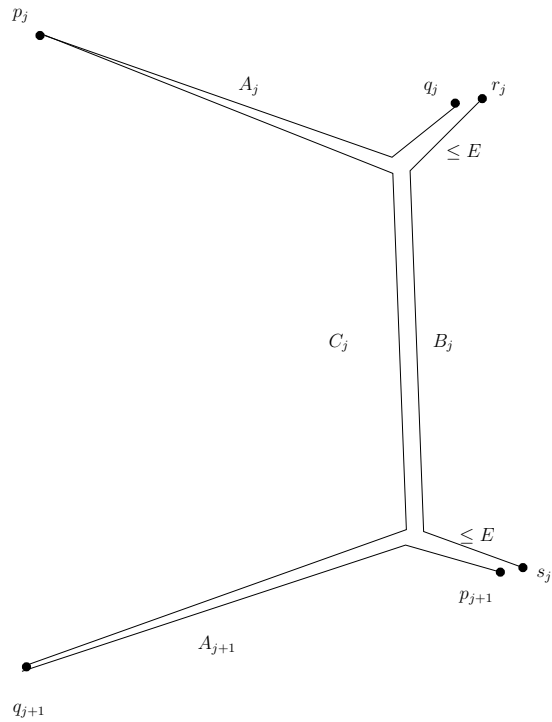
\begin{figure}[htbp]
\begin{center}
\scalebox{0.6}{\input{fig2.pstex_t}}
\caption{$A_j \dcap A_{j+1} = \emptyset$}
\label{figure2}
\end{center}
\end{figure}

For $j \ge 1$, we define
$$C_j=[p_j,q_{j+1}], \, D_j=[p_j,p_{j+1}].$$
The conditions (3),(4),(5) imply that
$d(C_j,B_j) \le 2E$ for all $j$.
It follows using (3) and (4) that 
$d(q_j,C_j) \le 4E, d(p_{j+1},C_j) \le 4E$.
This is because, for each $j$, the geodesic quadrilateral with the corners
$p_j,q_j,p_{j+1},q_{j+1}$ is $2 \delta$-thin
(see Figure \ref{figure1}, \ref{figure2}).

Therefore, 
for all $1 \le j \le i$, 
$|D_j \dcap D_{j+1}| \le 5E$.
It implies that for all $1 \le j \le i-1$,
$$|p_j-p_{j+2}| \ge \max(|p_j-p_{j+1}|, |p_{j+1}-p_{j+2}|) +990 E.$$
We checked the three points condition
for the constant $990E$ and the sequence  $\{p_j\}$, therefore  
by Proposition \ref{gromov.qi},  since $E>100 \delta$, we get that 
$|p_1-p_{i+1}| >0 $.

We are left to verify the condition (\ref{*}).
(1) follows from 
$|x-y| \le 2 \delta$.
Since $n_j \not=0$ and $|A_j|=\tr(a^{n_j})$, we have 
$|A_j| \ge \tr(a) \ge 1000 E$. We get (2).

Both $|A_j \dcap B_j|$ and 
$|B_{j-1} \dcap A_{j}|$
are at most $|\alpha \dcap \beta|$, 
therefore $\le D \le E$. We get (3) and (4).

To show (5), suppose not, i.e.,
$|A_j \dcap A_{j+1}| >E$ for some $j \ge 1$.
Set 
$$w_{j-1}=a^{n_1}b^{m_1} \cdots a^{n_{j-1}} b^{m_{j-1}}, \,
w_j=a^{n_1}b^{m_1} \cdots a^{n_{j}} b^{m_{j}}$$
(If $j=1$, then set $w_0=1$.)
Then $A_j$ is contained 
in $w_{j-1} (\alpha)$, which is 
a quasi-axis of $w_{j-1} f w_{j-1}^{-1}$, and 
also $A_{j+1}$ is 
contained in $w_{j} (\alpha)$, which is 
a quasi-axis of $w_{j} f w_{j}^{-1}$.
Let $\ell=[u,v] \subset A_{j+1}$ be a geodesic of length $E$
which is contained in $A_j \dcap A_{j+1}$.
For concreteness, suppose
the point $u$ is mapped toward
$v$ by (sufficiently big) positive powers of $w_{j-1} f w_{j-1}^{-1}$, 
and the point $v$ is mapped toward
$u$ by  (sufficiently big) positive powers of $w_{j} f w_{j}^{-1}$.
(Otherwise we take the inverse of 
the elements in the following.)
The direction of the actions makes sense since 
the set  $A_j \dcap A_{j+1}$ looks 
like a long narrow tube, which contains $\ell$.
Note that $\tr(w_{j} f w_{j}^{-1})
=\tr(w_{j-1} f w_{j-1}^{-1})=\tr(f)$
and $10 \tr(f) PKL  \le |\ell|$.
(See Figure \ref{figure4}.)
\begin{figure}[htbp]
\begin{center}
\scalebox{0.6}{\input{fig4.pstex_t}}
\caption{}
\label{figure4}
\end{center}
\end{figure}
We apply the following elements to the point $v$.
$$w_{j-1} f w_{j-1}^{-1} w_{j} f w_{j}^{-1}, \,
w_{j-1} f^2 w_{j-1}^{-1} w_{j} f^2 w_{j}^{-1}, 
\cdots,
w_{j-1} f^{PKL} w_{j-1}^{-1} w_{j} f^{PKL} w_{j}^{-1}.$$
Each of the elements moves $v$ at most $10 \delta$.
Set $v'= w_{j} f^{PKL} w_{j}^{-1}(v)$.
Then, each of the above elements also moves $v'$ 
at most $10 \delta$.
On the other hand
$|v-v'| = PKL \tr(f) \ge PKL/P \ge L$.
By the acylindricity, 
there must be
$1 \le I < J \le PKL$ such that 
$$w_{j-1} f^I w_{j-1}^{-1} w_{j} f^I w_{j}^{-1}
=w_{j-1} f^J w_{j-1}^{-1} w_{j} f^J w_{j}^{-1}$$
Since $w_{j-1}^{-1}w_j=a^{n_j}b^{m_j}$,
this implies that 
$[f^{I-J},a^{n_j}b^{m_j}]=1$.
Since $a=f^n$, 
it follows  that 
$[f^{I-J},b^{m_j}]=1$.
This is a contradiction (since $b=g$).
We showed (5).
The case (I) is completed.

For the case (II), we show $|y -w(y)|>0$. 
The argument is similar, and we omit the details.
So far, we have shown that $\langle g,f^n \rangle$ is 
free.

\noindent
{\bf Part 2}.
Now, we show that the embedding of $\langle g,f^n \rangle$
in $\Gamma$ by an orbit is quasi-isometric.
Our situation is same as in the proof of Proposition \ref{nielsen}.
Since the upper bound is trivial as before, we prove a 
lower bound.
Note that it is enough to get a desired uniform lower
bound for the case (O), (I) and (II) separately.
The case (O) is trivial, and the argument
is similar for (I) and (II).
We only discuss the case (I) in detail.

\begin{figure}[htbp]
\begin{center}
\scalebox{0.7}{\input{fig5.pstex_t}}
\caption{New points between $p_j,q_j$ ($n_j>0$)}
\label{figure5}
\end{center}
\end{figure}
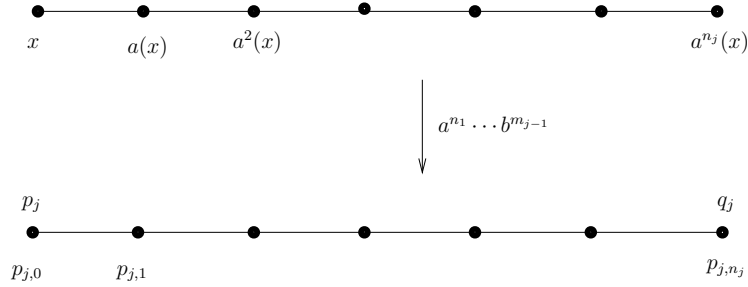

As in the proof of Proposition \ref{nielsen}, 
the reason why the above argument 
does not give a desired uniform lower bound in terms of $|w|$
is that $|p_{j+1}-p_j|$, the length of $D_j$, is unbounded when we vary $w$.
As before, we introduce interpolating points such that 
there is a uniform (for all $w$) upper bound on the distance 
between two consecutive points, and the three points 
condition is satisfied for $\varepsilon=E>100\delta$.
Between $p_j$ and $q_j$ for each $j$, we define points 
using the action of $f^n=a$ as follows:

\noindent
if $n_j>0$, define (see Figure \ref{figure5})
\begin{align*}
p_j=p_{j,0}&=a^{n_1}b^{m_1} \cdots a^{n_{j-1}} b^{m_{j-1}} (x), \\
p_{j,1} &= a^{n_1}b^{m_1} \cdots a^{n_{j-1}} b^{m_{j-1}} a(x),\\
& \cdots \\
p_{j,n_j} &=a^{n_1}b^{m_1} \cdots a^{n_{j-1}} b^{m_{j-1}}a^{n_j} (x)
=q_j.
\end{align*}
Note that for each $j$, the distance between any two 
consecutive points is $\tr(a)$.

If $n_j<0$, then define points similarly by 
$p_{j,k}= a^{n_1}b^{m_1} \cdots a^{n_{j-1}} b^{m_{j-1}} a^{-k}(x)$
for $0 \le k \le -n_j$.

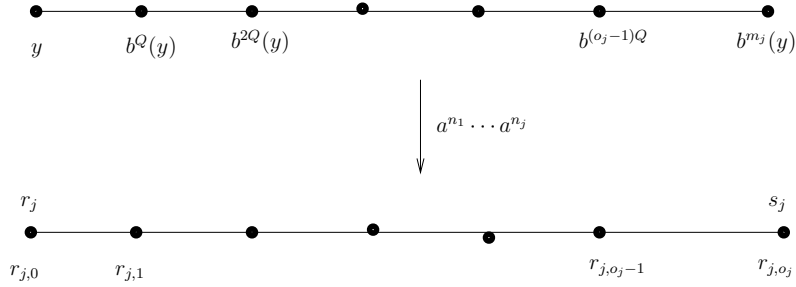
\begin{figure}[htbp]
\begin{center}
\scalebox{0.7}{\input{fig6.pstex_t}}
\caption{New points between $r_j,s_j$ ($m_j>0$)}
\label{figure6}
\end{center}
\end{figure}

Next, we also define points between $r_j$ and $s_j$.
In order to choose them with control on the distance between 
two consecutive new points, we fix an integer $Q \ge 1$ such that 
$$ (10E \le\,) \, \frac{\tr(a)}{100} \le \tr(b^Q) \le \frac{\tr(a)}{50}.$$
Such $Q$, which depends on $n$, exists since 
$\tr(b) \le \tr(f) \le \frac{\tr(a)}{100}$.
This is the only place where the condition 
$\tr(g) \le \tr(f)$ is used. (Remember $b=g$.)
Note that we use this $Q$ for all $j$.
To define points, we write (uniquely) for each $j$, 
$$m_j=o_jQ+l_j, $$
such that $o_j \in {\Bbb Z}$ and that 
$0 \le l_j <Q$ if $m_j \ge 0$ and
$-Q < l_j \le 0$ if $m_j <0$.
For each $j$, define points as follows between $r_j$ and $s_j$:

\noindent
if $m_j \ge 0$ and $o_j \ge 2$ then (see Figure \ref{figure6})
\begin{align*}
r_j=r_{j,0}&=a^{n_1}b^{m_1} \cdots a^{n_{j-1}} b^{m_{j-1}} a^{n_j}(y), \\
r_{j,1} &= a^{n_1}b^{m_1} \cdots a^{n_{j-1}} b^{m_{j-1}} a^{n_j}b^Q(y),\\
& \cdots \\
r_{j,o_j-1} &=a^{n_1}b^{m_1} \cdots a^{n_{j-1}} b^{m_{j-1}}a^{n_j} 
b^{(o_j-1)Q}(y) \\
r_{j,o_j} &=a^{n_1}b^{m_1} \cdots a^{n_{j-1}} b^{m_{j-1}}a^{n_j} 
b^{m_j}(y) 
=s_j.
\end{align*}

Here, if $o_j \not=0$, then 
the distance between any two consecutive points is $\tr(b^Q)$ 
except for the last pair of points, and 
it is between $\tr(b^Q)$ and $2\tr(b^Q)$ for the last pair.
If $o_j=1$, we do not produce any new points.
In this case, $\tr(b^Q) \le |r_j-s_j| \le 2\tr(b^Q)$.
Also, if $o_j=0$, then we do not produce any new points.
In that case, $|r_j -s_j| < \tr(b^Q)$.

If $m_j < 0$, we define points similarly, by
$r_{j,k}= a^{n_1}b^{m_1} \cdots a^{n_{j-1}} b^{m_{j-1}} a^{n_j}b^{-kQ}(y)$
for $0 \le k \le -o_j -1$ if $-o_j \ge 2$. Otherwise 
we do not produce new points.

In this way, we obtain a new set of points
$\{p_{j,k},r_{j,k}\}$ with the order which 
it naturally inherits from the order on the sequence 
$$p_1,q_1,r_1,s_1, \cdots, p_i,q_i,r_i,s_i.$$
The distance of any two consecutive points in this 
new sequence is at most $\tr(a)$.
But as a trade-off, the new sequence may not satisfy the three
points condition any more.
Therefore we modify the sequence by removing points.
First, remove all $r_j$.
($r_j$  is at most $2\delta$-close to $q_j$.)
Next, remove all $s_j$ except for the last point $s_i$. 
($s_j$  is at most $2\delta$-close to $p_{j+1}$.)
Finally, if $o_j=0$, then remove $q_j$.
We get a subsequence of points with order,
which we denote by $\{u_k\}$.

\begin{figure}[htbp]
\begin{center}
\scalebox{0.7}{\input{fig7.0.pstex_t}}
\caption{$q_j,r_j,s_j$ are removed if $o_j=0$}
\label{figure7.0}
\end{center}
\end{figure}
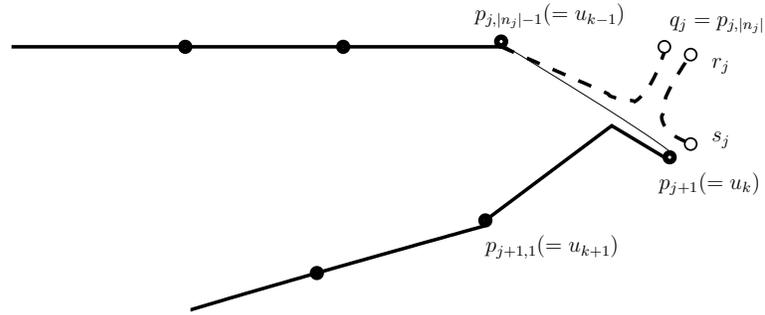

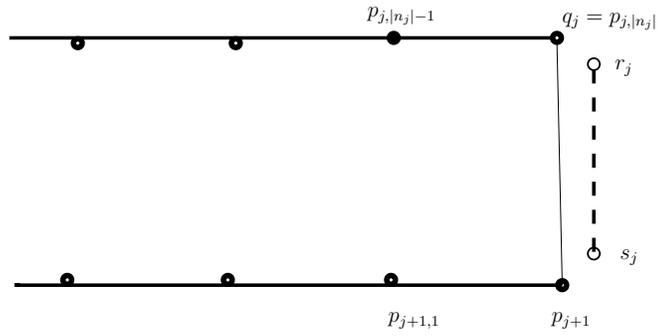
\begin{figure}[htbp]
\begin{center}
\scalebox{0.7}{\input{fig7.1.pstex_t}}
\caption{$r_j,s_j$ are removed if $o_j=1$}
\label{figure7.1}
\end{center}
\end{figure}

See Figure \ref{figure7.0},\ref{figure7.1},\ref{figure7.2}.
In those figures, two consecutive points
$u_k,u_{k+1}$ are joined by 
a solid line, where an interval by a thin solid 
line appears after we remove points while 
an interval by a thick solid line exists before 
we remove points. The intervals of dashed line do not 
exist because we removed points.
Removed points are described by white dots in the figures.

\begin{figure}[htbp]
\begin{center}
\scalebox{0.7}{\input{fig7.2.pstex_t}}
\caption{$r_j,s_j$ are removed if $o_j \ge 2$}
\label{figure7.2}
\end{center}
\end{figure}

Now, we argue that  the sequence $\{u_k\}$
satisfies the three points condition of Proposition 
\ref{gromov.qi}, and that there is a uniform
upper bound on the distance between any two 
consecutive points in the sequence.
In the following argument, it would help to 
keep in mind the following rough estimate of 
constants.

$$\delta < < E < < \tr(b^Q) < <\tr(a).$$

First, we remark that we 
removed points such that the distance of any two consecutive points
in $\{u_k\}$ is in one of the following intervals, in particular
there is a uniform upper bound:

\noindent
(i) between $\frac{4\tr(a)}{5}$ and $\frac{6\tr(a)}{5}$

\noindent
(ii) between $\frac{\tr(a)}{100}$ and $\frac{\tr(a)}{10}$.

Moreover,

\noindent
(iii) the sequence of intervals (between two consecutive points)
starts with 
$|n_1|$ intervals of length (i), followed by $|o_1|$ intervals 
of length (ii), 
followed by $|n_2|$ intervals of length (i), followed by 
$|o_2|$ intervals of length (ii), $\cdots$, 
$|n_i|$ intervals of length  (i), then $|o_i|$ intervals of length (ii)
at the end.
In particular, each interval has length 
$\ge \frac{\tr(a)}{100} \ge 10E \ge 1000$.

This is because removing $r_j$ and $s_j$ changes the 
distance between two consecutive points at most $4\delta$.
And removing $q_j$ happens only when $o_j=0$, therefore
$|r_j-s_j| \le \tr(b^Q) \le \frac{\tr(a)}{50}$.
In conclusion, the distance between any two consecutive 
points in the subsequence is 
at most $\tr(a)+4\delta + \frac{\tr(a)}{50} \le \frac{6}{5} \tr(a)$.

Next, we argue that the sequence  $\{u_k\}$
satisfies the three points condition of Proposition \ref{gromov.qi}
for the constant $E$.
This is by $\delta$-hyperbolicity, and essentially
by the same reason as the sequence $D_j$ satisfies the 
condition in the previous discussion(see Figure \ref{figure1}, \ref{figure2}).
If $o_j \ge 1$ (Figure \ref{figure7.1}, \ref{figure7.2}),
the three points (in other words, two consecutive intervals)
condition is nearly obvious since 
an interval of thin line is $10\delta$-neighborhood
of the corresponding thick dashed line interval.
The less obvious case is when $o_j=0$.
See Figure \ref{figure7.0}. We explain 
why the three points condition is satisfied
for $u_{k-1}=p_{j,|n_j|-1},u_k=p_{j+1},u_{k+1}=p_{j+1,1}$.
It suffices to show 
that $|[u_{k-1},u_k]\dcap [u_k,u_{k+1}]|$
is at most $\frac{1}{2}\min(|u_{k-1}-u_k|,|u_k-u_{k+1}|) -E$,
which is $\ge \frac{1}{2}\frac{4}{5}\tr(a)-E \ge \frac{39}{100}\tr(a)$.
First, $[u_{k-1},u_k]$ is contained
in the $10\delta$-neighborhood of 
$[u_{k-1},q_j]\cup[r_j,s_j]$ by 
$\delta$-hyperbolicity.
Therefore, $[u_{k-1},u_k]\dcap [u_k,u_{k+1}]$
is contained in the $10\delta$-neighborhood
of $([u_{k-1},q_j]\cup[r_j,s_j]) \dcap [u_k,u_{k+1}]$,
which is contained in $20\delta$-neighborhood of 
$([u_{k-1},q_j] \dcap [u_k,u_{k+1}])\cup [r_j,s_j]$.
But we know that 
$|[u_{k-1},q_j] \dcap [u_k,u_{k+1}]| \le E \le \frac{\tr(a)}{1000}$
and $|r_j-s_j| \le \frac{\tr(a)}{50}$.
The first inequality is by one of the conditions
we obtained, which is $|A_j \dcap A_{j+1}| \le E$
for all $j$.
Therefore, since $|r_j-q_j| \le 2 \delta$, 
the diameter of the union of those two sets
is smaller than $\frac{39}{100}\tr(a)$, which shows
our claim.
We verified the three points condition for $E$
for the sequence $\{u_k\}$, which starts at $x$ and 
ends at $w(x)$.

By Proposition \ref{gromov.qi} (for the first inequality below), we get
\begin{align*}
\lambda |w(x)-x| & \ge \sum_k |u_k-u_{k+1}| \ge 1000
\sum_{j=1}^i(|n_j|+ |o_j|) \ge \sum_{j=1}^i 500(|n_j|+ |o_j|+1)  \\
&\ge \frac{500}{Q} \sum_{j=1}^i(|n_j|+|m_j|)
= \frac{500}{Q}|w|,
\end{align*}
where $\lambda=(\frac{E}{100}-\delta)^{-1} \max_k|u_k-u_{k+1}|$.
The second inequality follows from the remark (iii) above.
We have the third inequality  because $|n_j| \ge 1$ for 
all $j$, and the fourth one because 
$Q(|o_j|+1) \ge |m_j|$ for all $j$.
Moreover, since $|u_k-u_{k+1}|$ is at most 
$2\tr(a)$, we have 
$\lambda \le (\frac{E}{100}-\delta)^{-1}2\tr(a)=\lambda_0$ 
for all $w$. 
Set $L=\frac{Q\lambda_0}{500}$, then 
we get $L|w(x)-x| \ge |w|$ for all $w$.
$L$ depends on $a,b$ but not on $w$ since 
so does $Q$.  
We completed the case (I). 
Note that we showed that the embedding is 
bi-Lipschitz for the point $x$ and the word metric
$|w|$ for the collection of words $w$ of the case (I).

The case (II) is similar, and we omit the details.
We finished the argument under the assumption that $D>0$.

Finally, assume that $D=0$.
The argument will be  essentially same as the case that 
$D>0$, so 
we discuss only the part which is different
(see the proof of Proposition \ref{nielsen}.
We discussed the case $D=0$ then the case $D>0$ as well).
Take the constants $K,L,P,E,N$ as before.
Set $a,b$ as before.
Now let $\ell$ be a geodesic which realizes
the distance between $\alpha,\beta$.
Let $m,x,y$ be the midpoint of $\ell$.
Then argue as before in the rest.
We define segments $A_j,B_j$ which satisfy 
the conditions (1)-(5). (In this case,
$q_j=r_j$ and $p_j=s_{j-1}$.)
This case is easier in the sense that 
since $\alpha \dcap \beta =\emptyset$,
we have $A_j \dcap B_j=\emptyset$, 
$B_j \dcap A_{j+1}=\emptyset$, 
therefore the union of those geodesics, which we need
to analyze in terms of metric, looks nearly 
like a tree with those segments as edges.
We omit the details.
\qed

Before we state a main theorem, we 
state a proposition which can be shown
similarly to Proposition \ref{nielsen2}.
The conclusion is weaker since $N_{\ref{nielsen3}}$ depends on $f,g$, but 
we do not require the conditions 1 and 2 in Proposition 
\ref{nielsen2} regarding  $\tr(f),\tr(g)$.

\begin{prop}\label{nielsen3}
Suppose $G$ acts on a $\delta$-hyperbolic graph $\Gamma$
acylindrically with constants $K(R),L(R)$.
Assume that $f,g \in G$  act hyperbolically
with quasi-axes $\alpha,\beta$.
Assume that $[f^s,g^t]\not=1$ for any $s,t \not=0$.

Then there exists a constant $N_{\ref{nielsen3}}=N>0$, 
which depends on $f,g$, such that if $n \ge N$, then 
$\langle g,f^n \rangle$ is free of rank two.
Moreover, the embedding 
of $\langle g,f^n \rangle$ in $\Gamma$ by an orbit is 
quasi-isometric. 
\end{prop}

\proof
The argument is very similar to the proof of Proposition \ref{nielsen2}, 
and easier since the difficult part was 
to obtain a uniform constant $N_{\ref{nielsen2}}$
for all $f,g$.
We only indicate where we need to modify the argument.

Choose constants $K,L,P$ in the same way.
Let $D=|\alpha \dcap \beta|$.
By Lemma \ref{axis}, $D < \infty$.
Define the constant $E$ in the same way.
Now fix a constant $N>0$ such that if $n \ge N$, then 
$$\tr(f^n) \ge 1000E, \text{ and } 
\tr(g) \le \frac{\tr(f^n)}{100}.$$
The constant $N$ depends on $f,g$ as opposed to 
Proposition \ref{nielsen2}, where we could 
choose $N_{\ref{nielsen2}}$ uniformly on $f,g$ 
to have those two inequalities for all $n \ge N_{\ref{nielsen2}}$
because of the conditions 1 and 2.
Note that the condition 2 (and condition 1) in Proposition \ref{nielsen2}
was used only to choose $N_{\ref{nielsen2}}$ uniformly on $f,g$ 
such that we have the first (and the second, respectively)
inequality in the above for all $n \ge N_{\ref{nielsen2}}$.

Let $n \ge N$ be a constant. Then we have the two
inequalities in the above.
Set $a=f^n,b=g$ as before. 
Then, using the first inequality, 
we can show that $\langle a,b \rangle$
is free exactly same as for Proposition \ref{nielsen2}.
What is essential in the argument is that 
$\tr(a)$ is much larger than $D$,
the $10 \delta$-overlap of the quasi-axes of $a,b$.
(Remember that $E \ge D$.)
We remark that  the condition 1 was irrelevant until this part
since it does not matter even if $\tr(g)=\tr(b)$ is much larger 
than $\tr(f)$ or $\tr(a)=\tr(f^n)$.

To show that the embedding is quasi-isometric for Proposition \ref{nielsen2},
we chose a constant $Q$
such that $\frac{\tr(a)}{100} \le \tr(b^Q) \le \frac{\tr(a)}{50}$.
This was possible since 
$\tr(b) \le \frac{\tr(a)}{100}$. (For this we used the condition 1.)
This is exactly the second inequality in the above, and 
we have chosen $N$ to have this inequality if $n \ge N$.
With this $Q$ we apply the same argument
for the rest.
We omit details.
\qed

\subsection{Upper bound on both exponents}\label{section.upperbound}
The following is the main theorem of Section \ref{section.freesubgroup}.
\begin{thm}\label{free}
Suppose $G$ acts acylindrically for constants
$K(R), L(R)$ on a $\delta$-hyperbolic graph $\Gamma$.
Then there exists a constant $M_{\ref{free}}$, which depends only
on $\delta$ and $K(20\delta),L(20\delta),L(200\delta)$
with the following property.

Suppose $a,b\in G$ act hyperbolically
with quasi-axes $\alpha,\beta$.
Assume for any $p,q \not=0$,
$[a^p,b^q] \not= 1$ in $G$.
Then for any $n,m \ge M$,
$\langle a^n,b^m \rangle$ is free of rank two.
Moreover, the embedding of $\langle a^n,b^m \rangle$ by an orbit 
in $\Gamma$ is quasi-isometric.
In particular, all non-trivial elements
in $\langle a^n,b^m \rangle$ are hyperbolic on $\Gamma$.
\end{thm}

%

\proof
Set $K=K(20\delta),L=L(20\delta)$.
Take constants $P_{\ref{min.tr}}$,
$N_{\ref{nielsen1}}$ and $N_{\ref{nielsen2}}$
by Lemma \ref{min.tr}, Proposition
\ref{nielsen1} and Proposition \ref{nielsen2}.
Set 
$$M=10KLP_{\ref{min.tr}}N_{\ref{nielsen1}}
+2000(\delta+1)P_{\ref{min.tr}}+N_{\ref{nielsen2}}.$$
$M$  depends only
on $\delta$ and $K(20\delta),L(20\delta),L(200\delta)$ and 
does not depend on $a,b$.
It suffices to show 

\noindent
{\bf Claim}.
$\langle a^n,b^m \rangle$ is free 
and the embedding to $\Gamma$ is quasi-isometric
if $n,m \ge M$.

To show this, set $D=|\alpha \dcap \beta|$.
Without loss of generality, we may assume
that $\tr(b) \le \tr(a)$.
By Lemma \ref{axis},
$D < 4PKL\tr(a) + 100\delta < \infty$.

\noindent
{\it Step 1}.
We may assume 
$\frac{\tr(b)}{\tr(a)} \le \frac{N_{\ref{nielsen1}}}{M}$.

\noindent
This is because otherwise we can show Claim
as follows.
Assume
$\frac{\tr(b)}{\tr(a)} > \frac{N_{\ref{nielsen1}}}{M}$.
Set $q=\frac{M}{N_{\ref{nielsen1}}}$,
and apply Proposition \ref{nielsen1} to $a,b$.
Then, if $n \ge N_{\ref{nielsen1}}$ and 
$m \ge q N_{\ref{nielsen1}}=M$, then 
$\langle a^n,b^m \rangle$ is free and 
the embedding is quasi-isometric.
Since $M > N_{\ref{nielsen1}}$, we get Claim.

Note that it follows that $10KLP \tr(b) \le \tr(a)$
since 
$\frac{N_{\ref{nielsen1}}}{M} < \frac{1}{10KLP}$.

\noindent
{\it Step 2}.
We may assume
$M\tr(b) \le D+100(\delta +1)$.

\noindent
This is because otherwise we get Claim
as follows. Assume the contrary.
It immediately follows that 
both $\tr(a^n),\tr(b^m)$ are $> D+100(\delta+1)$.
(We use $\tr(b) \le \tr(a)$.)
Then by Proposition \ref{nielsen}, 
we get Claim. 

Note that since 
$10KL+2000(\delta +1) \le M\tr(b)$ by the way 
we chose $M$, it follows that 
$10KL+1000 \delta \le  D$.

\noindent
{\it Step 3}.
We have $D \le 2 \tr(a)$
(assuming the inequalities in Step 1 and 2).

\noindent
To argue by contradiction, assume $D > 2 \tr(a)$.
Then we get a contradiction  using the same idea 
as for Lemma \ref{axis} concerning 
the action of commutators $[b^i,a]$.
Since $D \ge 1000 \delta$ by Step 2,
the set $\alpha \dcap \beta$ looks
like a narrow tube. Therefore, it makes sense
to talk about the direction of the action by $a$ and $b$
along this tube, and furthermore, the direction of $a$
coincides the direction of one of $b$ or $b^{-1}$.
In the following, we assume that the actions by $a,b$ 
have the same direction along $\alpha \dcap \beta$, otherwise,
we consider $b^{-1}$ instead of $b$.

Let $\ell=[p,p'] \subset \alpha$ be the longest 
segment contained in $\alpha \cap N_{2\delta}(\beta)$
such that $a$ moves $p$ toward $p'$,
i.e., $a(p) \in [p,p']$. We know $|\ell| \ge D -100 \delta$.
Since $D \ge 1000 \delta$ by Step 2,
it follows $|\ell| \ge \frac{9}{10}D$.

We claim that for all $1 \le i \le PKL$,
we have $d(p,[b^i,a](p)) \le 20 \delta$.
If $\delta=0$ then this is obvious. 
Indeed, first of all, $\alpha$ and $\beta$ coincide in $\ell$
in this case.
Also, if we apply $a,b^i,a^{-1}$ then $b^{-i}$ in this 
order to $p$, the point moves in $\ell$.
This is because since 
$10KLP \tr(b) \le \tr(a)$ by Step 1 and 
$D > 2 \tr(a)$,
we have $\tr(a)+\tr( b^i) \le \frac{11}{10} \tr(a) 
\le \frac{11}{20} D \le |\ell|$ 
for all $1 \le i \le PKL$.
Now it is trivial that $p=[b^i,a](p)$.
If $\delta >0$, when 
we apply $a,b^i,a^{-1}$ then $b^{-i}$ to $p$, the point moves 
in the $10\delta$-neighborhood of $\ell$, and 
we get $d(p,[b^i,a](p)) \le 20 \delta$
for all $1 \le i \le PKL$. (see Fig. \ref{comm}
with the roles of $a$ and $b$ exchanged.)

Let $q \in \ell$ be the point with $L = d(p,q)$.
(Note that $L \le \frac{D}{10}$ by Step 2.)
By the same reason as for $p$, 
we have $d(q,[b^i,a](q)) \le 20 \delta$
for all $1 \le i \le PKL$.
Now by the acylindricity, since $d(p,q)=L=L(20\delta)$, 
there must be $1 \le I < J \le PKL$ such that 
$[b^I,a] =[b^J,a]$,
therefore $[b^{I-J},a]=1$, 
which gives a contradiction.
This is the end of Step 3.

But if $D \le 2 \tr(a)$,
we can apply Proposition \ref{nielsen2} to $a,b$. 
Therefore, if $n \ge N_{\ref{nielsen2}}$, then 
$\langle b,a^n\rangle $ is free and 
the embedding is quasi-isometric.
Since $M > N_{\ref{nielsen2}}$,
we have shown the Claim. Note that $\langle b^m, a^n \rangle$
is a subgroup of $\langle b, a^n\rangle $.
\qed

\begin{remark}
It is more difficult to deal with 
the normal subgroup generated by $a^n,b^m$, or 
even just by $a^n$ (see Question 11 \cite{Iv2}).
See the work of Delzant \cite{De}.

\end{remark}

\begin{remark}
Theorem \ref{free} is regarding two elements, but 
one can ask if there exists  a constant $M$ such that 
if $a,b,c \in G$ are hyperbolic elements with certain 
condition (for example, pairwise independence), 
then $\langle a^{\ell},b^m,c^n\rangle$
is free for any $\ell,m,n \ge M$.
We remark that the rank of the free subgroup
may not be three.
Take two hyperbolic elements $a,b \in G$ which 
satisfy the commutator assumption in Theorem \ref{free}
(i.e., independent).
For any $M>0$, set $c=a^M b a^{-M}$.
Then, the pairs $a,c$ and $b,c$ are also
independent, but $\langle a^M,b^M,c^M \rangle$ is 
equal to $\langle a^M,b^M \rangle$.
\end{remark}

\subsection{Hyperbolic isometries without quasi-axes}\label{no.axes}
So far, we have been discussing hyperbolic isometries
with quasi-axes. 
The existence of quasi-axes, which are geodesics by our 
definition, is a restriction but indeed not really necessary
for our arguments since we can use certain quasi-geodesics
and modify the original arguments. 
We discuss this issue in this section. 
Readers may skip this section since we do not use this
for our main application to mapping class groups in Section \ref{section.mcg}.

A path $\alpha$ parametrized by 
the arc-length is called a {\it $(K,\ep)$-quasi-geodesic}
for $0<K \le 1$ and $0 \le \ep$ 
if for all $t,s$ we have 
$K|t-s|-\ep \le d(\alpha(t),\alpha(s))$.

The following fact is elementary
(see \cite{Fu} for details).
\begin{fact}[quasi-geodesic axis]\label{qg.axis}
If $a$ is a hyperbolic isometry of 
a $\delta$-hyperbolic graph $\Gamma$, there exists a $(K,\ep)$-quasi-geodesic $\alpha$
for some $K,\ep$ such that 
\begin{enumerate}
\item
$a^n(\alpha)$ and $\alpha$ are in the $30 \delta$-neighborhood
of each other for any $n$. (Namely, $\alpha$ is almost invariant by $a$.)
\item
Let $p,q \in \alpha$.
Then the subpath of $\alpha$ between
$p,q$ and a geodesic $[p,q]$ are in the $10 \delta$-neighborhood 
of each other.
\end{enumerate}
We call such path $\alpha$ as a {\it quasi-geodesic axis} of $a$
in this paper. To be precise, we should use the term
{\it quasi-geodesic quasi-axis}, but we make it shorter.
\end{fact}
One can easily show from (1) and (2) that any two quasi-geodesic axes of
$a$ are in the $30 \delta$-neighborhood of each other.
Note that (2) is concerning only the path, but not the element $a$.
Also, the quasi-geodesic constants of $\alpha$ are not important
for our purpose. What is useful for us is (2).

For example, Lemma \ref{min.tr} gives a uniform positive 
lower bound of $\tr(a)$ for all hyperbolic isometry $a \in G$
with a quasi-axis if the action of $G$ is acylindrical, 
but, indeed the assumption on the existence 
of quasi-axes is redundant.
The proof of Lemma \ref{min.tr} easily generalizes
by using quasi-geodesic axes instead of quasi-axes
(see \cite{Fu} for the precise argument), and we obtain
the following.
\begin{lemma}\label{min.tr.noaxes}
Suppose $G$ acts on a $\delta$-hyperbolic graph
$\Gamma$.
If the action is acylindrical with constants
$K(R)$ and $L(R)$, then there exists
an integer $P \ge 1$ such that for any element
$a \in G$ which acts hyperbolically on $\Gamma$, 
we have $\tr(a^P) \ge 1$.
The constant $P$ depends only on $\delta$
and $K(200\delta)$.
\end{lemma}
Of course, this constant $P$ is maybe larger than
$P_{\ref{min.tr}}$.
The existence of such $P$, but not the actual number, 
is essential for our argument.

We restate Theorem \ref{free}  in the following form
for a potential application.
The only difference is that we do not assume that 
there are quasi-axes for $a$ and $b$.
\begin{thm}\label{free'}
Suppose $G$ acts acylindrically for constants
$K(R), L(R)$ on a $\delta$-hyperbolic graph $\Gamma$.
Then there exists a constant $M'_{\ref{free'}}$, which depends only
on $\delta$ and $K(20\delta),L(20\delta),L(200\delta)$,
with the following property.

Suppose $a,b\in G$ act hyperbolically.
Assume for any $p,q \not=0$,
$[a^p,b^q] \not= 1$ in $G$.
Then for any $n,m \ge M'$,
$\langle a^n,b^m \rangle$ is free of rank two.
Moreover, the embedding of $\langle a^n,b^m \rangle$ by an orbit 
in $\Gamma$ is quasi-isometric.
In particular, all non-trivial elements
in $\langle a^n,b^m \rangle$ are hyperbolic on $\Gamma$.
\end{thm}

\proof
The proof is very similar to the one for 
Theorem \ref{free}.
Basically, we use quasi-geodesic axes instead of 
quasi-axes for hyperbolic isometries.
The proof of Theorem  \ref{free}
relies on Proposition \ref{nielsen}, Lemma \ref{min.tr}, Lemma \ref{axis},
Proposition \ref{nielsen1} and Proposition \ref{nielsen2}.
We have already generalized Lemma \ref{min.tr}
to Lemma \ref{min.tr.noaxes}.
We modify the statement and the proof
of each of the other ones, which we only outline here.

As for Proposition \ref{nielsen},
replace quasi-axes $\alpha, \beta$ by
quasi-geodesic axes $\alpha, \beta$. Accordingly, replace
all $\alpha \dcap \beta$ by $\alpha \cap_{1000\delta} \beta$
in the statement. Then the original proof works with minor modification
using the properties (1) and (2) in Fact \ref{qg.axis}
of the quasi-geodesic axes.
 
In Lemma \ref{axis},
replace quasi-axes $\alpha, \beta$ for $a,b$
by quasi-geodesic axes in the assumption.
We also replace the inequality in the conclusion by
$$|\alpha \cap _{1000 \delta} \beta| < 4P_{\ref{min.tr.noaxes}}
K(20\delta) L(20\delta)
\max(\tr(a),\tr(b)) + 10000 \delta.$$
Then, the proof is same after an appropriate modification
regarding constants.

As for Proposition \ref{nielsen1}, Proposition \ref{nielsen2},
replace quasi-axes $\alpha, \beta$ by quasi-geodesic axes.
Use $\alpha \cap_{1000\delta} \beta$ instead of 
$\alpha \dcap \beta$ in the statement.
Then the original proof works with minor modification. 
Having done them all, we modify the proof
of Theorem \ref{free} to fit our setting.
Of course, we always use the constant $P_{\ref{min.tr.noaxes}}$ 
instead of $P_{\ref{min.tr}}$ .
The constant $M'_{\ref{free'}}$ is maybe
larger than $M_{\ref{free}}$.
We omit details.
\qed

\begin{remark}
Proposition \ref{nielsen3} also holds
if we drop the assumption on the existence 
of quasi-axes $\alpha, \beta$. 
The argument is also very similar to the original one.
\end{remark}

\section{Application to mapping class group}\label{section.mcg}
We discuss mapping class groups in this 
section. We apply results from Section 
\ref{section.freesubgroup}
to pseudo-Anosov elements.
Theorem \ref{mod} and \ref{mod2} are
main results of the paper. 
\subsection{Uniform estimate}
We apply Theorem \ref{free} to 
the mapping class group, $\Mod(S)$, of a 
compact orientable surface $S$.
Let $\C(S)$ be the curve graph of $S$
(see for example \cite{Iv}, \cite{MM} for the definition).
Masur-Minsky \cite{MM} showed that 
$\C(S)$ is $\delta$-hyperbolic and 
an element $a \in \Mod(S)$
is pseudo-Anosov if and only if it 
acts as a hyperbolic isometry on $\C(S)$, and
moreover (\cite{Bo}) there always exists a quasi-axis.
Bowditch \cite{Bo} showed that the action is 
acylindrical.

For a subgroup $G < \Mod(S)$, Farb-Mosher \cite{farb.mosher}
introduced the notion of {\it convex-cocompact}.
It has been shown (\cite{Ha}, \cite{KeLe}) that $G$ is convex-cocompact iff
for a point $c \in \C(S)$, the map from $G$ to $\C(S)$
sending $g$ to $g(c)$, namely the embedding by 
an orbit,  is quasi-isometric.

The following is an immediate consequence of Theorem \ref{free}.
Apply it to the action of $\Mod(S)$ to $\C(S)$.
Two pseudo-Anosov elements
$a,b$ are called {\it independent}
if  $[a^n,b^m]\not=1$ for any $n,m \not=0$ (cf.\cite{Iv}).

\begin{thm}\label{mod}
Let $S$ be a compact orientable surface, and
$\Mod(S)$ its mapping class group.
Then there exists a constant $M(S)$ with the following property.
Suppose $a,b \in \Mod(S)$ are pseudo-Anosov elements
such that $[a^n,b^m]\not=1$ for any $n,m \not=0$.
Then for any $n,m \ge M$, 
$\langle a^n,b^m \rangle $ is free of rank two, and 
convex-cocompact. In particular 
all non-trivial elements in $\langle a^n,b^m \rangle $
are pseudo-Anosov.
\end{thm}



\subsection{Non-uniform estimate and example}

Let $a,b \in \Mod(S)$ be two independent pseudo-Anosov 
elements.
It would be interesting to know for which $(n,m)$, 
$\langle a^n,b^m \rangle$ is 
free of rank two, and convex-cocompact.
The following theorem says that 
it is the case except for finitely many $(n,m)$.
We do not know if the number of the exceptional pairs
is bounded.

\begin{thm}\label{mod2}
Let $S$ be a compact orientable 
surface and $a,b$ two independent
pseudo-Anosov elements.
Then there exists $N_{\ref{mod2}}=N$, which depends on $a,b$, such that for any
$n \ge N$, both 
$\langle a,b^n \rangle$ and  $\langle b,a^n \rangle$
are free of rank two, and convex-cocompact.
In particular, 
$\langle a^n,b^m \rangle$ is 
free of rank two, and convex-cocompact
if $|n|+|m| \ge 2N$ and $nm\not=0$.
\end{thm}

\proof
Apply Proposition \ref{nielsen3}
to $a,b$ for the action on $\C(S)$.
The constant $N_{\ref{nielsen3}}$ will do.
\qed

The constant $N_{\ref{mod2}}$ must depend on $a,b$
as the following example shows.
\begin{example}\label{example}
Let $S$ be a compact orientable surface
which is not a sphere with less than four punctures
or a torus. 
If  $n>0$ is sufficiently large, then there exist
two independent pseudo-Anosov elements $f,g \in \Mod(S)$
such that 
$\langle g,f^n \rangle $ is not free.

To see this, take $f,a \in \Mod(S)$ such that 
$f$ is pseudo-Anosov, $a$ is non-trivial torsion
and that $\langle f,a \rangle$ is not virtually cyclic.
To find such $a,f$, first take a
non-trivial torsion element $a \in \Mod(S)$
such that there is a non-trivial and non-peripheral 
simple closed curve $\sigma$ on $S$ which is not 
homotopic to $a(\sigma)$. One can find such $a$ easily.
Then one can find a desired $f$.
For example, take any pseudo-Anosov element $h$
on $S$.
Let $d$ be a Dehn-twist along $\sigma$.
Set $f=d^m h d^{-m}$. We choose a sufficiently large $m>0$ later.
It is clear that $f$ is pseudo-Anosov, and 
the two laminations which are invariant by $f$,
which we regard as a set of two points, ${\rm fix}(f)$, in the 
boundary of the Teichmuller space of $S$, must
be moved by $a$ (i.e. ${\rm fix}(f) \cap a({\rm fix}(f))=\emptyset$)
 if $m$ is sufficiently large. For such $m$, 
it follows by a standard
argument that $\langle f,a \rangle$ is not virtually
cyclic (cf. \cite{Iv1}).

Now, for sufficiently large $n$, $f^na$ is pseudo-Anosov, 
and independent from $f$.
One can show this using the curve graph of $S$, 
$\C(S)$, which is $\delta$-hyperbolic.
If necessary, replace $f$ by some power of it 
in advance, and we may assume that $f$ leaves
a geodesic $\gamma$ in $\C(S)$ invariant.
By our assumption $\gamma \dcap a(\gamma)$
is bounded.
For each $n>0$, one can find a line 
which is invariant by $f^na$ 
using a piece of $\gamma$, a fundamental 
domain for the action of $f^n$, and the action of $a$.
Then, for sufficiently large $n$, using
$\delta$-hyperbolic geometry of $\C(S)$, one can 
show that the line is indeed a quasi-geodesic, 
therefore $f^na$ is pseudo-Anosov.
Moreover, for sufficiently large $n$, 
the quasi-geodesic has two points at infinity of $\C(S)$
which are disjoint from the two points for $\gamma$.
It implies that $f^na$ and $f$ are independent.
Set $g=f^na$.
Then $\langle g,f^n \rangle$ is not free since 
it contains the torsion element $a$.
\end{example}

\end{document}

%% file: fig0.pstex_t
\begin{picture}(0,0)%
\includegraphics{fig0.pstex}%
\end{picture}%
\setlength{\unitlength}{4144sp}%
\begingroup\makeatletter\ifx\SetFigFontNFSS\undefined%
\gdef\SetFigFontNFSS#1#2#3#4#5{%
  \reset@font\fontsize{#1}{#2pt}%
  \fontfamily{#3}\fontseries{#4}\fontshape{#5}%
  \selectfont}%
\fi\endgroup%
\begin{picture}(7812,2404)(616,-2249)
\put(4276,-556){\makebox(0,0)[lb]{\smash{{\SetFigFontNFSS{12}{14.4}{\rmdefault}{\mddefault}{\updefault}{\color[rgb]{0,0,0}$p$}%
}}}}
\put(2656,-1771){\makebox(0,0)[lb]{\smash{{\SetFigFontNFSS{12}{14.4}{\rmdefault}{\mddefault}{\updefault}{\color[rgb]{0,0,0}$\beta$}%
}}}}
\put(2746,-421){\makebox(0,0)[lb]{\smash{{\SetFigFontNFSS{12}{14.4}{\rmdefault}{\mddefault}{\updefault}{\color[rgb]{0,0,0}$\alpha$}%
}}}}
\put(1036,-106){\makebox(0,0)[lb]{\smash{{\SetFigFontNFSS{12}{14.4}{\rmdefault}{\mddefault}{\updefault}{\color[rgb]{0,0,0}$A^{-1}(p)$}%
}}}}
\put(631,-2176){\makebox(0,0)[lb]{\smash{{\SetFigFontNFSS{12}{14.4}{\rmdefault}{\mddefault}{\updefault}{\color[rgb]{0,0,0}$B^{-1}(p)$}%
}}}}
\put(7966,-2131){\makebox(0,0)[lb]{\smash{{\SetFigFontNFSS{12}{14.4}{\rmdefault}{\mddefault}{\updefault}{\color[rgb]{0,0,0}$B(p)$}%
}}}}
\put(7246,-16){\makebox(0,0)[lb]{\smash{{\SetFigFontNFSS{12}{14.4}{\rmdefault}{\mddefault}{\updefault}{\color[rgb]{0,0,0}$A(p)$}%
}}}}
\put(4951,-691){\makebox(0,0)[lb]{\smash{{\SetFigFontNFSS{12}{14.4}{\rmdefault}{\mddefault}{\updefault}{\color[rgb]{0,0,0}$\ell$}%
}}}}
\end{picture}%

%% file: ueg.fig3.pstex_t
\begin{picture}(0,0)%
\includegraphics{ueg.fig3.pstex}%
\end{picture}%
\setlength{\unitlength}{4144sp}%
\begingroup\makeatletter\ifx\SetFigFontNFSS\undefined%
\gdef\SetFigFontNFSS#1#2#3#4#5{%
  \reset@font\fontsize{#1}{#2pt}%
  \fontfamily{#3}\fontseries{#4}\fontshape{#5}%
  \selectfont}%
\fi\endgroup%
\begin{picture}(8259,2416)(799,-2240)
\put(8281,-1726){\makebox(0,0)[lb]{\smash{{\SetFigFontNFSS{12}{14.4}{\rmdefault}{\mddefault}{\updefault}{\color[rgb]{0,0,0}$\beta$}%
}}}}
\put(6301,-646){\makebox(0,0)[lb]{\smash{{\SetFigFontNFSS{12}{14.4}{\rmdefault}{\mddefault}{\updefault}{\color[rgb]{0,0,0}$a^{4PKL}(p)$}%
}}}}
\put(5716,-691){\makebox(0,0)[lb]{\smash{{\SetFigFontNFSS{12}{14.4}{\rmdefault}{\mddefault}{\updefault}{\color[rgb]{0,0,0}$\ell$}%
}}}}
\put(2521,-646){\makebox(0,0)[lb]{\smash{{\SetFigFontNFSS{12}{14.4}{\rmdefault}{\mddefault}{\updefault}{\color[rgb]{0,0,0}$p$}%
}}}}
\put(3151,-646){\makebox(0,0)[lb]{\smash{{\SetFigFontNFSS{12}{14.4}{\rmdefault}{\mddefault}{\updefault}{\color[rgb]{0,0,0}$x=a^{PKL}(p)$}%
}}}}
\put(4276,-646){\makebox(0,0)[lb]{\smash{{\SetFigFontNFSS{12}{14.4}{\rmdefault}{\mddefault}{\updefault}{\color[rgb]{0,0,0}$y=a^{2PKL}(p)$}%
}}}}
\put(8416,-2176){\makebox(0,0)[lb]{\smash{{\SetFigFontNFSS{12}{14.4}{\rmdefault}{\mddefault}{\updefault}{\color[rgb]{0,0,0}$b (\text{or }b^{-1})$}%
}}}}
\put(8281,-421){\makebox(0,0)[lb]{\smash{{\SetFigFontNFSS{12}{14.4}{\rmdefault}{\mddefault}{\updefault}{\color[rgb]{0,0,0}$a$}%
}}}}
\put(8056, 29){\makebox(0,0)[lb]{\smash{{\SetFigFontNFSS{12}{14.4}{\rmdefault}{\mddefault}{\updefault}{\color[rgb]{0,0,0}$\alpha$}%
}}}}
\end{picture}%

%% file: comm.pstex_t
\begin{picture}(0,0)%
\includegraphics{comm.pstex}%
\end{picture}%
\setlength{\unitlength}{4144sp}%
\begingroup\makeatletter\ifx\SetFigFontNFSS\undefined%
\gdef\SetFigFontNFSS#1#2#3#4#5{%
  \reset@font\fontsize{#1}{#2pt}%
  \fontfamily{#3}\fontseries{#4}\fontshape{#5}%
  \selectfont}%
\fi\endgroup%
\begin{picture}(7269,976)(1609,-1475)
\put(8326,-646){\makebox(0,0)[lb]{\smash{{\SetFigFontNFSS{12}{14.4}{\rmdefault}{\mddefault}{\updefault}{\color[rgb]{0,0,0}$\alpha$}%
}}}}
\put(8461,-1276){\makebox(0,0)[lb]{\smash{{\SetFigFontNFSS{12}{14.4}{\rmdefault}{\mddefault}{\updefault}{\color[rgb]{0,0,0}$\beta$}%
}}}}
\put(3871,-1366){\makebox(0,0)[lb]{\smash{{\SetFigFontNFSS{12}{14.4}{\rmdefault}{\mddefault}{\updefault}{\color[rgb]{0,0,0}$a^{-i}ba^i(x)$}%
}}}}
\put(2566,-691){\makebox(0,0)[lb]{\smash{{\SetFigFontNFSS{12}{14.4}{\rmdefault}{\mddefault}{\updefault}{\color[rgb]{0,0,0}$x$}%
}}}}
\put(1891,-1231){\makebox(0,0)[lb]{\smash{{\SetFigFontNFSS{12}{14.4}{\rmdefault}{\mddefault}{\updefault}{\color[rgb]{0,0,0}$b^{-1}a^{-i}ba^i(x)$}%
}}}}
\put(7291,-1411){\makebox(0,0)[lb]{\smash{{\SetFigFontNFSS{12}{14.4}{\rmdefault}{\mddefault}{\updefault}{\color[rgb]{0,0,0}$ba^i(x)$}%
}}}}
\put(5851,-646){\makebox(0,0)[lb]{\smash{{\SetFigFontNFSS{12}{14.4}{\rmdefault}{\mddefault}{\updefault}{\color[rgb]{0,0,0}$a^i(x)$}%
}}}}
\end{picture}%

%% file: fig1.pstex_t
\begin{picture}(0,0)%
\includegraphics{fig1.pstex}%
\end{picture}%
\setlength{\unitlength}{4144sp}%
\begingroup\makeatletter\ifx\SetFigFontNFSS\undefined%
\gdef\SetFigFontNFSS#1#2#3#4#5{%
  \reset@font\fontsize{#1}{#2pt}%
  \fontfamily{#3}\fontseries{#4}\fontshape{#5}%
  \selectfont}%
\fi\endgroup%
\begin{picture}(5025,3484)(391,-3464)
\put(2701,-1186){\makebox(0,0)[lb]{\smash{{\SetFigFontNFSS{12}{14.4}{\rmdefault}{\mddefault}{\updefault}{\color[rgb]{0,0,0}$A_j$}%
}}}}
\put(406,-3391){\makebox(0,0)[lb]{\smash{{\SetFigFontNFSS{12}{14.4}{\rmdefault}{\mddefault}{\updefault}{\color[rgb]{0,0,0}$q_{j+1}$}%
}}}}
\put(586,-151){\makebox(0,0)[lb]{\smash{{\SetFigFontNFSS{12}{14.4}{\rmdefault}{\mddefault}{\updefault}{\color[rgb]{0,0,0}$p_j$}%
}}}}
\put(2116,-2536){\makebox(0,0)[lb]{\smash{{\SetFigFontNFSS{12}{14.4}{\rmdefault}{\mddefault}{\updefault}{\color[rgb]{0,0,0}$A_{j+1}$}%
}}}}
\put(4681,-871){\makebox(0,0)[lb]{\smash{{\SetFigFontNFSS{12}{14.4}{\rmdefault}{\mddefault}{\updefault}{\color[rgb]{0,0,0}$q_j$}%
}}}}
\put(4681,-2491){\makebox(0,0)[lb]{\smash{{\SetFigFontNFSS{12}{14.4}{\rmdefault}{\mddefault}{\updefault}{\color[rgb]{0,0,0}$p_{j+1}$}%
}}}}
\put(5401,-871){\makebox(0,0)[lb]{\smash{{\SetFigFontNFSS{12}{14.4}{\rmdefault}{\mddefault}{\updefault}{\color[rgb]{0,0,0}$r_j$}%
}}}}
\put(5311,-2446){\makebox(0,0)[lb]{\smash{{\SetFigFontNFSS{12}{14.4}{\rmdefault}{\mddefault}{\updefault}{\color[rgb]{0,0,0}$s_j$}%
}}}}
\put(4726,-1681){\makebox(0,0)[lb]{\smash{{\SetFigFontNFSS{12}{14.4}{\rmdefault}{\mddefault}{\updefault}{\color[rgb]{0,0,0}$B_j$}%
}}}}
\put(5041,-1996){\makebox(0,0)[lb]{\smash{{\SetFigFontNFSS{12}{14.4}{\rmdefault}{\mddefault}{\updefault}{\color[rgb]{0,0,0}$\le E$}%
}}}}
\put(3646,-1501){\makebox(0,0)[lb]{\smash{{\SetFigFontNFSS{12}{14.4}{\rmdefault}{\mddefault}{\updefault}{\color[rgb]{0,0,0}$\le 2E$}%
}}}}
\put(2341,-1591){\makebox(0,0)[lb]{\smash{{\SetFigFontNFSS{12}{14.4}{\rmdefault}{\mddefault}{\updefault}{\color[rgb]{0,0,0}$C_j$}%
}}}}
\put(4996,-1456){\makebox(0,0)[lb]{\smash{{\SetFigFontNFSS{12}{14.4}{\rmdefault}{\mddefault}{\updefault}{\color[rgb]{0,0,0}$\le E$}%
}}}}
\end{picture}%

%% file: fig2.pstex_t
\begin{picture}(0,0)%
\includegraphics{fig2.pstex}%
\end{picture}%
\setlength{\unitlength}{4144sp}%
\begingroup\makeatletter\ifx\SetFigFontNFSS\undefined%
\gdef\SetFigFontNFSS#1#2#3#4#5{%
  \reset@font\fontsize{#1}{#2pt}%
  \fontfamily{#3}\fontseries{#4}\fontshape{#5}%
  \selectfont}%
\fi\endgroup%
\begin{picture}(5250,7174)(571,-7154)
\put(586,-151){\makebox(0,0)[lb]{\smash{{\SetFigFontNFSS{12}{14.4}{\rmdefault}{\mddefault}{\updefault}{\color[rgb]{0,0,0}$p_j$}%
}}}}
\put(4681,-871){\makebox(0,0)[lb]{\smash{{\SetFigFontNFSS{12}{14.4}{\rmdefault}{\mddefault}{\updefault}{\color[rgb]{0,0,0}$q_j$}%
}}}}
\put(586,-7081){\makebox(0,0)[lb]{\smash{{\SetFigFontNFSS{12}{14.4}{\rmdefault}{\mddefault}{\updefault}{\color[rgb]{0,0,0}$q_{j+1}$}%
}}}}
\put(2431,-6451){\makebox(0,0)[lb]{\smash{{\SetFigFontNFSS{12}{14.4}{\rmdefault}{\mddefault}{\updefault}{\color[rgb]{0,0,0}$A_{j+1}$}%
}}}}
\put(3736,-3436){\makebox(0,0)[lb]{\smash{{\SetFigFontNFSS{12}{14.4}{\rmdefault}{\mddefault}{\updefault}{\color[rgb]{0,0,0}$C_j$}%
}}}}
\put(2836,-871){\makebox(0,0)[lb]{\smash{{\SetFigFontNFSS{12}{14.4}{\rmdefault}{\mddefault}{\updefault}{\color[rgb]{0,0,0}$A_j$}%
}}}}
\put(4771,-3436){\makebox(0,0)[lb]{\smash{{\SetFigFontNFSS{12}{14.4}{\rmdefault}{\mddefault}{\updefault}{\color[rgb]{0,0,0}$B_j$}%
}}}}
\put(4906,-1546){\makebox(0,0)[lb]{\smash{{\SetFigFontNFSS{12}{14.4}{\rmdefault}{\mddefault}{\updefault}{\color[rgb]{0,0,0}$\le E$}%
}}}}
\put(5041,-5371){\makebox(0,0)[lb]{\smash{{\SetFigFontNFSS{12}{14.4}{\rmdefault}{\mddefault}{\updefault}{\color[rgb]{0,0,0}$\le E$}%
}}}}
\put(5806,-5731){\makebox(0,0)[lb]{\smash{{\SetFigFontNFSS{12}{14.4}{\rmdefault}{\mddefault}{\updefault}{\color[rgb]{0,0,0}$s_j$}%
}}}}
\put(5041,-5956){\makebox(0,0)[lb]{\smash{{\SetFigFontNFSS{12}{14.4}{\rmdefault}{\mddefault}{\updefault}{\color[rgb]{0,0,0}$p_{j+1}$}%
}}}}
\put(5356,-871){\makebox(0,0)[lb]{\smash{{\SetFigFontNFSS{12}{14.4}{\rmdefault}{\mddefault}{\updefault}{\color[rgb]{0,0,0}$r_j$}%
}}}}
\end{picture}%

%% file: fig4.pstex_t
\begin{picture}(0,0)%
\includegraphics{fig4.pstex}%
\end{picture}%
\setlength{\unitlength}{4144sp}%
\begingroup\makeatletter\ifx\SetFigFontNFSS\undefined%
\gdef\SetFigFontNFSS#1#2#3#4#5{%
  \reset@font\fontsize{#1}{#2pt}%
  \fontfamily{#3}\fontseries{#4}\fontshape{#5}%
  \selectfont}%
\fi\endgroup%
\begin{picture}(7584,2364)(979,-2413)
\put(6481,-1636){\makebox(0,0)[lb]{\smash{{\SetFigFontNFSS{12}{14.4}{\rmdefault}{\mddefault}{\updefault}{\color[rgb]{0,0,0}$v$}%
}}}}
\put(4771,-1636){\makebox(0,0)[lb]{\smash{{\SetFigFontNFSS{12}{14.4}{\rmdefault}{\mddefault}{\updefault}{\color[rgb]{0,0,0}$w_jf^kw_j^{-1}(v)$}%
}}}}
\put(1801,-601){\makebox(0,0)[lb]{\smash{{\SetFigFontNFSS{12}{14.4}{\rmdefault}{\mddefault}{\updefault}{\color[rgb]{0,0,0}$w_{j-1}fw_{j-1}^{-1}$}%
}}}}
\put(5626,-691){\makebox(0,0)[lb]{\smash{{\SetFigFontNFSS{12}{14.4}{\rmdefault}{\mddefault}{\updefault}{\color[rgb]{0,0,0}$w_{j-1}f^kw_{j-1}w_jf^kw_j^{-1}(v)$}%
}}}}
\put(8461,-331){\makebox(0,0)[lb]{\smash{{\SetFigFontNFSS{12}{14.4}{\rmdefault}{\mddefault}{\updefault}{\color[rgb]{0,0,0}$A_j \subset w_{j-1}(\alpha)$}%
}}}}
\put(8371,-2086){\makebox(0,0)[lb]{\smash{{\SetFigFontNFSS{12}{14.4}{\rmdefault}{\mddefault}{\updefault}{\color[rgb]{0,0,0}$A_{j+1} \subset w_j(\alpha)$}%
}}}}
\put(2431,-1591){\makebox(0,0)[lb]{\smash{{\SetFigFontNFSS{12}{14.4}{\rmdefault}{\mddefault}{\updefault}{\color[rgb]{0,0,0}$u$}%
}}}}
\put(1576,-2086){\makebox(0,0)[lb]{\smash{{\SetFigFontNFSS{12}{14.4}{\rmdefault}{\mddefault}{\updefault}{\color[rgb]{0,0,0}$w_jfw_j^{-1}$}%
}}}}
\end{picture}%

%% file: fig5.pstex_t
\begin{picture}(0,0)%
\includegraphics{fig5.pstex}%
\end{picture}%
\setlength{\unitlength}{4144sp}%
\begingroup\makeatletter\ifx\SetFigFontNFSS\undefined%
\gdef\SetFigFontNFSS#1#2#3#4#5{%
  \reset@font\fontsize{#1}{#2pt}%
  \fontfamily{#3}\fontseries{#4}\fontshape{#5}%
  \selectfont}%
\fi\endgroup%
\begin{picture}(6151,2407)(526,-2564)
\put(1531,-601){\makebox(0,0)[lb]{\smash{{\SetFigFontNFSS{12}{14.4}{\familydefault}{\mddefault}{\updefault}{\color[rgb]{0,0,0}$a(x)$}%
}}}}
\put(676,-556){\makebox(0,0)[lb]{\smash{{\SetFigFontNFSS{12}{14.4}{\familydefault}{\mddefault}{\updefault}{\color[rgb]{0,0,0}$x$}%
}}}}
\put(2431,-556){\makebox(0,0)[lb]{\smash{{\SetFigFontNFSS{12}{14.4}{\familydefault}{\mddefault}{\updefault}{\color[rgb]{0,0,0}$a^2(x)$}%
}}}}
\put(6346,-556){\makebox(0,0)[lb]{\smash{{\SetFigFontNFSS{12}{14.4}{\familydefault}{\mddefault}{\updefault}{\color[rgb]{0,0,0}$a^{n_j}(x)$}%
}}}}
\put(631,-1906){\makebox(0,0)[lb]{\smash{{\SetFigFontNFSS{12}{14.4}{\familydefault}{\mddefault}{\updefault}{\color[rgb]{0,0,0}$p_j$}%
}}}}
\put(6571,-1906){\makebox(0,0)[lb]{\smash{{\SetFigFontNFSS{12}{14.4}{\familydefault}{\mddefault}{\updefault}{\color[rgb]{0,0,0}$q_j$}%
}}}}
\put(541,-2491){\makebox(0,0)[lb]{\smash{{\SetFigFontNFSS{12}{14.4}{\familydefault}{\mddefault}{\updefault}{\color[rgb]{0,0,0}$p_{j,0}$}%
}}}}
\put(1441,-2491){\makebox(0,0)[lb]{\smash{{\SetFigFontNFSS{12}{14.4}{\familydefault}{\mddefault}{\updefault}{\color[rgb]{0,0,0}$p_{j,1}$}%
}}}}
\put(6481,-2446){\makebox(0,0)[lb]{\smash{{\SetFigFontNFSS{12}{14.4}{\familydefault}{\mddefault}{\updefault}{\color[rgb]{0,0,0}$p_{j,n_j}$}%
}}}}
\put(4186,-1276){\makebox(0,0)[lb]{\smash{{\SetFigFontNFSS{12}{14.4}{\familydefault}{\mddefault}{\updefault}{\color[rgb]{0,0,0}$a^{n_1}\cdots b^{m_{j-1}}$}%
}}}}
\end{picture}%

%% file: fig6.pstex_t
\begin{picture}(0,0)%
\includegraphics{fig6.pstex}%
\end{picture}%
\setlength{\unitlength}{4144sp}%
\begingroup\makeatletter\ifx\SetFigFontNFSS\undefined%
\gdef\SetFigFontNFSS#1#2#3#4#5{%
  \reset@font\fontsize{#1}{#2pt}%
  \fontfamily{#3}\fontseries{#4}\fontshape{#5}%
  \selectfont}%
\fi\endgroup%
\begin{picture}(6691,2407)(526,-2564)
\put(721,-601){\makebox(0,0)[lb]{\smash{{\SetFigFontNFSS{12}{14.4}{\familydefault}{\mddefault}{\updefault}{\color[rgb]{0,0,0}$y$}%
}}}}
\put(1531,-601){\makebox(0,0)[lb]{\smash{{\SetFigFontNFSS{12}{14.4}{\familydefault}{\mddefault}{\updefault}{\color[rgb]{0,0,0}$b^Q(y)$}%
}}}}
\put(2431,-556){\makebox(0,0)[lb]{\smash{{\SetFigFontNFSS{12}{14.4}{\familydefault}{\mddefault}{\updefault}{\color[rgb]{0,0,0}$b^{2Q}(y)$}%
}}}}
\put(6751,-556){\makebox(0,0)[lb]{\smash{{\SetFigFontNFSS{12}{14.4}{\familydefault}{\mddefault}{\updefault}{\color[rgb]{0,0,0}$b^{m_j}(y)$}%
}}}}
\put(5401,-556){\makebox(0,0)[lb]{\smash{{\SetFigFontNFSS{12}{14.4}{\familydefault}{\mddefault}{\updefault}{\color[rgb]{0,0,0}$b^{(o_j-1)Q}$}%
}}}}
\put(4186,-1276){\makebox(0,0)[lb]{\smash{{\SetFigFontNFSS{12}{14.4}{\familydefault}{\mddefault}{\updefault}{\color[rgb]{0,0,0}$a^{n_1}\cdots a^{n_{j}}$}%
}}}}
\put(631,-1906){\makebox(0,0)[lb]{\smash{{\SetFigFontNFSS{12}{14.4}{\familydefault}{\mddefault}{\updefault}{\color[rgb]{0,0,0}$r_j$}%
}}}}
\put(7021,-1906){\makebox(0,0)[lb]{\smash{{\SetFigFontNFSS{12}{14.4}{\familydefault}{\mddefault}{\updefault}{\color[rgb]{0,0,0}$s_j$}%
}}}}
\put(541,-2491){\makebox(0,0)[lb]{\smash{{\SetFigFontNFSS{12}{14.4}{\familydefault}{\mddefault}{\updefault}{\color[rgb]{0,0,0}$r_{j,0}$}%
}}}}
\put(1441,-2491){\makebox(0,0)[lb]{\smash{{\SetFigFontNFSS{12}{14.4}{\familydefault}{\mddefault}{\updefault}{\color[rgb]{0,0,0}$r_{j,1}$}%
}}}}
\put(6931,-2446){\makebox(0,0)[lb]{\smash{{\SetFigFontNFSS{12}{14.4}{\familydefault}{\mddefault}{\updefault}{\color[rgb]{0,0,0}$r_{j,o_j}$}%
}}}}
\put(5491,-2446){\makebox(0,0)[lb]{\smash{{\SetFigFontNFSS{12}{14.4}{\familydefault}{\mddefault}{\updefault}{\color[rgb]{0,0,0}$r_{j,o_j-1}$}%
}}}}
\end{picture}%

%% file: fig7.0.pstex_t
\begin{picture}(0,0)%
\includegraphics{fig7.0.pstex}%
\end{picture}%
\setlength{\unitlength}{4144sp}%
\begingroup\makeatletter\ifx\SetFigFontNFSS\undefined%
\gdef\SetFigFontNFSS#1#2#3#4#5{%
  \reset@font\fontsize{#1}{#2pt}%
  \fontfamily{#3}\fontseries{#4}\fontshape{#5}%
  \selectfont}%
\fi\endgroup%
\begin{picture}(6033,2679)(598,-2209)
\put(6616,-106){\makebox(0,0)[lb]{\smash{{\SetFigFontNFSS{12}{14.4}{\familydefault}{\mddefault}{\updefault}{\color[rgb]{0,0,0}$r_j$}%
}}}}
\put(6616,-736){\makebox(0,0)[lb]{\smash{{\SetFigFontNFSS{12}{14.4}{\familydefault}{\mddefault}{\updefault}{\color[rgb]{0,0,0}$s_j$}%
}}}}
\put(6256,254){\makebox(0,0)[lb]{\smash{{\SetFigFontNFSS{12}{14.4}{\familydefault}{\mddefault}{\updefault}{\color[rgb]{0,0,0}$q_j=p_{j,|n_j|}$}%
}}}}
\put(4681,-1681){\makebox(0,0)[lb]{\smash{{\SetFigFontNFSS{12}{14.4}{\familydefault}{\mddefault}{\updefault}{\color[rgb]{0,0,0}$p_{j+1,1}(=u_{k+1})$}%
}}}}
\put(6166,-1141){\makebox(0,0)[lb]{\smash{{\SetFigFontNFSS{12}{14.4}{\familydefault}{\mddefault}{\updefault}{\color[rgb]{0,0,0}$p_{j+1}(=u_k)$}%
}}}}
\put(4591,299){\makebox(0,0)[lb]{\smash{{\SetFigFontNFSS{12}{14.4}{\familydefault}{\mddefault}{\updefault}{\color[rgb]{0,0,0}$p_{j,|n_j|-1}(=u_{k-1})$}%
}}}}
\end{picture}%

%% file: fig7.1.pstex_t
\begin{picture}(0,0)%
\includegraphics{fig7.1.pstex}%
\end{picture}%
\setlength{\unitlength}{4144sp}%
\begingroup\makeatletter\ifx\SetFigFontNFSS\undefined%
\gdef\SetFigFontNFSS#1#2#3#4#5{%
  \reset@font\fontsize{#1}{#2pt}%
  \fontfamily{#3}\fontseries{#4}\fontshape{#5}%
  \selectfont}%
\fi\endgroup%
\begin{picture}(5268,2854)(1498,-2384)
\put(4591,299){\makebox(0,0)[lb]{\smash{{\SetFigFontNFSS{12}{14.4}{\familydefault}{\mddefault}{\updefault}{\color[rgb]{0,0,0}$p_{j,|n_j|-1}$}%
}}}}
\put(6256,254){\makebox(0,0)[lb]{\smash{{\SetFigFontNFSS{12}{14.4}{\familydefault}{\mddefault}{\updefault}{\color[rgb]{0,0,0}$q_j=p_{j,|n_j|}$}%
}}}}
\put(6166,-2311){\makebox(0,0)[lb]{\smash{{\SetFigFontNFSS{12}{14.4}{\familydefault}{\mddefault}{\updefault}{\color[rgb]{0,0,0}$p_{j+1}$}%
}}}}
\put(4771,-2311){\makebox(0,0)[lb]{\smash{{\SetFigFontNFSS{12}{14.4}{\familydefault}{\mddefault}{\updefault}{\color[rgb]{0,0,0}$p_{j+1,1}$}%
}}}}
\put(6751,-1771){\makebox(0,0)[lb]{\smash{{\SetFigFontNFSS{12}{14.4}{\familydefault}{\mddefault}{\updefault}{\color[rgb]{0,0,0}$s_j$}%
}}}}
\put(6706,-151){\makebox(0,0)[lb]{\smash{{\SetFigFontNFSS{12}{14.4}{\familydefault}{\mddefault}{\updefault}{\color[rgb]{0,0,0}$r_j$}%
}}}}
\end{picture}%

%% file: fig7.2.pstex_t
\begin{picture}(0,0)%
\includegraphics{fig7.2.pstex}%
\end{picture}%
\setlength{\unitlength}{4144sp}%
\begingroup\makeatletter\ifx\SetFigFontNFSS\undefined%
\gdef\SetFigFontNFSS#1#2#3#4#5{%
  \reset@font\fontsize{#1}{#2pt}%
  \fontfamily{#3}\fontseries{#4}\fontshape{#5}%
  \selectfont}%
\fi\endgroup%
\begin{picture}(5178,4969)(1498,-4499)
\put(4906,299){\makebox(0,0)[lb]{\smash{{\SetFigFontNFSS{12}{14.4}{\familydefault}{\mddefault}{\updefault}{\color[rgb]{0,0,0}$p_{j,|n_j|-1}$}%
}}}}
\put(6256,299){\makebox(0,0)[lb]{\smash{{\SetFigFontNFSS{12}{14.4}{\familydefault}{\mddefault}{\updefault}{\color[rgb]{0,0,0}$q_j=p_{j,|n_j|}$}%
}}}}
\put(6616,-61){\makebox(0,0)[lb]{\smash{{\SetFigFontNFSS{12}{14.4}{\familydefault}{\mddefault}{\updefault}{\color[rgb]{0,0,0}$r_j$}%
}}}}
\put(6661,-3661){\makebox(0,0)[lb]{\smash{{\SetFigFontNFSS{12}{14.4}{\familydefault}{\mddefault}{\updefault}{\color[rgb]{0,0,0}$s_j$}%
}}}}
\put(4996,-4426){\makebox(0,0)[lb]{\smash{{\SetFigFontNFSS{12}{14.4}{\familydefault}{\mddefault}{\updefault}{\color[rgb]{0,0,0}$p_{j+1,1}$}%
}}}}
\put(6301,-4336){\makebox(0,0)[lb]{\smash{{\SetFigFontNFSS{12}{14.4}{\familydefault}{\mddefault}{\updefault}{\color[rgb]{0,0,0}$p_{j+1}$}%
}}}}
\put(6616,-1051){\makebox(0,0)[lb]{\smash{{\SetFigFontNFSS{12}{14.4}{\familydefault}{\mddefault}{\updefault}{\color[rgb]{0,0,0}$r_{j,1}$}%
}}}}
\put(6571,-2311){\makebox(0,0)[lb]{\smash{{\SetFigFontNFSS{12}{14.4}{\familydefault}{\mddefault}{\updefault}{\color[rgb]{0,0,0}$r_{j,o_j-1}$}%
}}}}
\end{picture}%